\documentclass[11pt]{amsart}

\usepackage[T1]{fontenc}
\usepackage{lmodern}
\usepackage{amsmath,amssymb,amsthm,mathtools,mathrsfs}
\usepackage{microtype}
\usepackage[margin=1.08in]{geometry}
\usepackage[hidelinks]{hyperref}
\usepackage{aliascnt}
\hypersetup{pdftitle={Zeros of Katz\textquotesingle s p-adic Eisenstein series E2},pdfauthor={Abdellah Sebbar and Khalil Besrour}}
\usepackage[nameinlink,capitalize]{cleveref}
\usepackage{enumitem}

\numberwithin{equation}{section}

\theoremstyle{plain}
\newtheorem{theorem}{Theorem}[section]
\newaliascnt{proposition}{theorem}
\newtheorem{proposition}[proposition]{Proposition}
\aliascntresetthe{proposition}
\newaliascnt{lemma}{theorem}
\newtheorem{lemma}[lemma]{Lemma}
\aliascntresetthe{lemma}
\newaliascnt{corollary}{theorem}
\newtheorem{corollary}[corollary]{Corollary}
\aliascntresetthe{corollary}
\newtheorem*{theoremA}{Theorem A}
\newtheorem*{theoremB}{Theorem B}
\newtheorem*{theoremC}{Theorem C}
\newtheorem*{theoremD}{Theorem D}
\theoremstyle{definition}
\newaliascnt{definition}{theorem}
\newtheorem{definition}[definition]{Definition}
\aliascntresetthe{definition}
\newaliascnt{remark}{theorem}
\newtheorem{remark}[remark]{Remark}
\aliascntresetthe{remark}
\newaliascnt{example}{theorem}

\aliascntresetthe{example}
\newaliascnt{question}{theorem}
\newtheorem{question}[question]{Question}
\aliascntresetthe{question}
\crefname{theorem}{Theorem}{Theorems}
\crefname{proposition}{Proposition}{Propositions}
\crefname{lemma}{Lemma}{Lemmas}
\crefname{corollary}{Corollary}{Corollaries}
\crefname{definition}{Definition}{Definitions}
\crefname{remark}{Remark}{Remarks}
\crefname{example}{Example}{Examples}
\crefname{question}{Question}{Questions}

\newcommand{\Q}{\mathbf{Q}}
\newcommand{\Qbar}{\overline{\mathbf Q}}
\newcommand{\Z}{\mathbf{Z}}
\newcommand{\C}{\mathbf{C}}
\newcommand{\Qp}{\mathbf{Q}_p}
\newcommand{\Zp}{\mathbf{Z}_p}
\newcommand{\Cp}{\mathbf{C}_p}
\newcommand{\Fp}{\mathbf{F}_p}
\newcommand{\Fbarp}{\overline{\mathbf{F}}_p}

\newcommand{\Hdg}{\mathcal H}
\newcommand{\Bpoly}{\mathcal B}
\newcommand{\thetaop}{\theta}

\newcommand{\dR}{\mathrm{dR}}

\begin{document}

\title[Zeros of Katz\textquotesingle s Eisenstein series]{Zeros of Katz's $p$-adic Eisenstein Series $E_2$:\\
Hypergeometric congruences and isogenies}
\author{Abdellah Sebbar}
\author{Khalil Besrour}
\address{Department of Mathematics and Statistics, University of Ottawa, Ottawa, Ontario K1N 6N5, Canada}
\email{asebbar@uottawa.ca}
\email{kbesr067@uottawa.ca}

\subjclass{Primary 11F33, 11F85; Secondary 11G07, 33C05.}
\keywords{$p$-adic modular forms, Eisenstein series, supersingular polynomial, Dwork congruences, Legendre family, unit-root splitting}

\begin{abstract}
For every prime $p\geq5$, the zeros of Katz's $p$-adic
Eisenstein series $E_2$ on elliptic curves with good ordinary
reduction have exactly $\lfloor p/12\rfloor+1$ distinct
$j$-invariants. We prove that all these invariants lie in finite
unramified extensions of $\Qp$. Their reductions are given by
a classical formula involving the supersingular polynomial,
and each nonelliptic residue class contains a unique simple
zero. For the Legendre family, hypergeometric congruences
determine the polynomial of the lifted invariants modulo
every power of $p$. We also study isogenies between the zeros.
Each $p$-power isogeny class contains at most one zero
$j$-invariant, and any two distinct isogenous zero invariants
are algebraic over $\Q$. At a fixed cyclic degree, the
nonelliptic isogenies give points in a finite divisor over
$\Q$ on the corresponding modular curve. Its degree bounds
their number and their fields of definition, and excludes
several small degrees. The only zeros with complex
multiplication have $j=0$ or $1728$. Whether a non-CM zero
can have an algebraic $j$-invariant remains open.
\end{abstract}
\maketitle
\enlargethispage{2pt}

\section{Introduction}

The normalized weight-two Eisenstein series
\[
 E_2(\tau)=1-24\sum_{n\geq1}\sigma_1(n)q^n,
 \qquad q=e^{2\pi i\tau},
\]
where $\sigma_1(n)=\sum_{d\mid n}d$, is quasimodular over
$\C$. Its Fourier expansion is also a $p$-adic modular form
in the sense of Serre \cite{SerrePadic}. Katz gives a geometric
interpretation of this form on elliptic curves with ordinary
reduction, using the unit-root splitting of de Rham cohomology
\cite{Katz1973,Katz1976}. This interpretation allows one to
study its zeros and the effect of isogenies on them.

In this paper, we determine these zeros and their fields of
definition. We then study isogenies between the corresponding
elliptic curves. The geometric criterion for vanishing is
particularly simple. For a normalized equation
$y^2=4x^3-g_2x-g_3$, put $\omega=dx/y$ and $\eta=x\,dx/y$.
Katz's $E_2$ vanishes precisely when the class of $\eta$
generates the unit-root line. We use this criterion together
with congruences for classical Eisenstein series to count the
zeros, and with V\'elu's formula to restrict the isogenies
between them.

Fix a prime $p\geq5$. We write $\Qp$ for the field of
$p$-adic numbers, $\Zp$ for its ring of integers, and $\Cp$
for the completion of an algebraic closure of $\Qp$.
The field with $p$ elements is denoted by $\Fp$, and its
algebraic closure by $\Fbarp$. A Katz zero will mean an
elliptic curve with good ordinary reduction on which
Katz's $E_2$ vanishes. Vanishing is unchanged when the
invariant differential is multiplied by a nonzero scalar,
so it is a condition on the $j$-invariant. We call the
values $j=0$ and $1728$ elliptic, and all other values
nonelliptic.

The reductions of the zeros can be expressed in terms of
the supersingular polynomial. Write
\[
 p-1=12m+4\delta_p+6\varepsilon_p,
 \qquad m=\left\lfloor\frac p{12}\right\rfloor,
\]
where
\[
 \delta_p=
 \begin{cases}
 0,&p\equiv1\pmod3,\\
 1,&p\equiv2\pmod3,
 \end{cases}
 \qquad
 \varepsilon_p=
 \begin{cases}
 0,&p\equiv1\pmod4,\\
 1,&p\equiv3\pmod4.
 \end{cases}
\]
Thus $\delta_p=1$ exactly when $j=0$ is supersingular, and $\varepsilon_p=1$ exactly when $j=1728$ is supersingular.  Let
\[
 S_p(X)=\prod_{j_s}(X-j_s)\in\Fp[X]
\]
be the monic polynomial whose roots are the distinct
supersingular $j$-invariants in $\Fbarp$, and put
\[
 Z_p(X)=-12S_p'(X)
 +8\delta_p\frac{S_p(X)}{X}
 +6\varepsilon_p\frac{S_p(X)}{X-1728}.
\]
Each quotient is a polynomial when its coefficient is
nonzero. This formula is classical. Kaneko and Zagier
obtain it, in their notation $\widetilde E_{p+1}(j)$, in
the proof of their Atkin congruence
\cite[\S6, pp.~112--113]{KanekoZagier1998}. Here $E_k$
denotes the normalized Eisenstein series of weight $k$.
The coprimality of the reductions of $E_{p-1}$ and $E_{p+1}$
is the statement used in Serre's account of Swinnerton-Dyer's
argument \cite[\S1.4, Theorem~5 and Corollary~1]{SerreCong}.
Katz also proves it by identifying the relevant Cartier
coefficients with the reductions of $E_{p-1}$ and $E_{p+1}/12$
\cite[Theorem~3.1]{KatzZannier2021}. A further factorization
is given in \cite[Lemma~2 and Table~1]{Bianchini2026}.

We show that the roots of $Z_p$ are simple. Hensel's lemma
then gives one zero of the Katz form above each root.
The uniqueness holds in the entire residue disc, including
points defined over ramified extensions. Together with
the ordinary elliptic zeros, this gives the following result.

\begin{theoremA}
The roots of $Z_p$ are ordinary, nonelliptic, and simple.  Each root has a unique simple lift to a zero of Katz's $p$-adic $E_2$, and a root of degree $e$ over $\Fp$ lifts over the unramified extension of $\Qp$ of exact degree $e$.  If $d=\deg Z_p$, the lifted $j$-invariants are the roots of a monic polynomial
\[
 \mathcal Z_p(X)\in\Zp[X],\qquad
 \mathcal Z_p(X)\equiv Z_p(X)\pmod p,
\]
whose discriminant is a $p$-adic unit. Hence the finite
$\Zp$-scheme defined by $\mathcal Z_p$ is \'{e}tale.
Including the ordinary elliptic classes, the zeros have
\[
 \left\lfloor\frac p{12}\right\rfloor+1
\]
distinct $j$-invariants.  When $p\equiv1\pmod{12}$, the nonelliptic reductions are the critical points of $S_p$.
\end{theoremA}

The factorization of $Z_p$ over $\Fp$ therefore determines
the Frobenius orbits of the nonelliptic zeros. The resulting
finite \'{e}tale algebra over $\Zp$ is determined by its
reduction. Its distinguished $j$-coordinate has characteristic
polynomial $\mathcal Z_p$.

To describe the lifted polynomial more explicitly, we pass
to the Legendre family $y^2=x(x-1)(x-\lambda)$, where
$\lambda\ne0,1$. Write
\[
 F(\lambda)={}_2F_1\!\left(\frac12,\frac12;1;\lambda\right),
 \qquad \mathcal D=(1-2\lambda)+6\lambda(1-\lambda)\frac{d}{d\lambda},
\]
and let
$\Hdg_p(\lambda)=\sum_{k=0}^{(p-1)/2}\binom{(p-1)/2}{k}^2\lambda^k$
be the Deuring polynomial over $\Fp$. Its roots are the
supersingular Legendre parameters. The operator $\mathcal D$
occurs in the classical identity $E_2=F\mathcal DF$.
We use its reduction modulo $p$ to compare the zeros in the
Legendre parameter with those in the $j$-coordinate.
The complex arithmetic consequences of the same identity
are studied separately in \cite{SebbarComplex}.

\begin{theoremB}
The normalized Katz form on the ordinary Legendre line reduces to
\[
 \frac{\overline E_2}{\overline C}
 =\frac{\mathcal D\Hdg_p}{\Hdg_p},\qquad C=F^2.
\]
The nonelliptic zero divisor is the pullback of the divisor of $Z_p$. Above an ordinary elliptic invariant, its multiplicity is two at each point over $j=0$ and one at each point over $j=1728$.
\end{theoremB}

The full factorization is given in
Theorem~\ref{thm:mainFactorization}. The logarithmic-derivative
congruence of Beukers and Vlasenko extends this comparison
modulo $p^s$ for every $s\ge1$. Put
\[
 F_s(T)=\sum_{n=0}^{p^s-1}\frac{(1/2)_n^2}{(n!)^2}T^n,
 \quad \Bpoly_{p,s}=\mathcal D F_s,
 \quad N(T)=256(T^2-T+1)^3,\quad D(T)=T^2(1-T)^2.
\]

Here $(a)_0=1$ and $(a)_n=a(a+1)\cdots(a+n-1)$ for $n\ge1$
denote the rising factorials. The map to the $j$-line is
$j=N(T)/D(T)$.

\begin{theoremC}
The quotient $\Bpoly_{p,s}/F_s$ agrees with the normalized
Katz form modulo $p^s$ on the ordinary part of the Legendre
family. Each nonelliptic residue zero has a unique root of $\Bpoly_{p,s}$ in its disc, congruent modulo $p^s$ to the Katz zero. If $G_{p,s}$ is the monic factor selecting these discs and $d=\deg Z_p$, then
\[
 G_{p,s}(T)\equiv256^{-d}D(T)^d\mathcal Z_p\bigl(N(T)/D(T)\bigr)\pmod{p^s}.
\]
The coefficients in this identity determine
$\mathcal Z_p\bmod p^s$ uniquely.
\end{theoremC}

Proposition~\ref{prop:KatzDworkComparison} proves the analytic
comparison, including the ordinary cusp discs.
Theorem~\ref{thm:direct} gives the relation between the two
polynomials and an explicit inverse relation between their
coefficients. When $d=0$, the nonelliptic zero polynomials
and the corresponding factors are empty products, equal to $1$.

Our second aim is to study isogenies between Katz zeros.
V\'elu's normalized formula expresses the change of $E_2$
under an isogeny as a sum of the $x$-coordinates of its
nonzero kernel points \cite[\S2.4]{Kohel1996}. At a fixed
cyclic degree $N$, the vanishing of this sum gives a finite
divisor over $\Q$ on the modular curve $X_0(N)$, which
parametrizes cyclic subgroups of order $N$. Let $g(X_0(N))$
denote its genus. We obtain the following restrictions.
Here CM means complex multiplication.

\begin{theoremD}
The only CM Katz zeros have $j=0$ or $1728$, with the corresponding ordinary condition. Each $p$-power isogeny class contains at most one Katz-zero $j$-invariant, and no nonelliptic zero admits a cyclic $p$-power isogeny to another zero.

No nonelliptic Katz zero admits a cyclic isogeny to another zero of degree $2,3,4,5,7,9,13$, or $25$. More generally, for each prime $\ell$, there are at most $g(X_0(\ell))$ unordered pairs of nonelliptic zero $j$-invariants joined by a cyclic degree-$\ell$ isogeny. For each such pair,
\[
 [\Q(j(E),j(E')):\Q]\leq2g(X_0(\ell)).
\]
Any two distinct isogenous Katz-zero $j$-invariants are algebraic over $\Q$.
\end{theoremD}

The restriction at $p$-power degree follows from the
unramified fields of definition and uniqueness of the lifted
zeros. A cyclic $p$-power isogeny between two nonelliptic
zeros would identify an iterate of the canonical quotient
with a Frobenius conjugate. Repeating the quotient would
then give a cyclic self-isogeny. This forces complex
multiplication, contrary to the classification of CM zeros.

For arbitrary degree $N>1$, the same sum vanishes, away
from elliptic points, precisely at the zeros of
\[
 d\log\frac{\Delta(\tau)}{\Delta(N\tau)}.
\]
Here $\Delta=(E_4^3-E_6^2)/1728$ is the modular
discriminant. The degree of the zero divisor of this
differential is expressed in terms of the genus and cusp
orders of $X_0(N)$. It bounds the number of isogenies,
the degree over $\Q$ of the field generated by the two
$j$-invariants, and the least common multiple of their
residue degrees over $\Fp$.

The distinction between algebraicity over $\Qp$ and
algebraicity over $\Q$ is essential here. Every zero invariant
is algebraic over $\Qp$, but it remains open whether a
non-CM zero invariant can be algebraic over $\Q$. We explain
in the final section how this question is related to the
$p$-adic period conjecture of Ancona and
Fr\u{a}\c{t}il\u{a} \cite{AnconaFratila2025}.
Their Conjecture~8.3 implies nonvanishing on each non-CM
algebraic curve to which it applies. Their Remark~7.7
verifies the additional unramified condition on the motive
in this case. This implication remains conditional.

Section~\ref{sec:Katz} establishes the unit-root criterion.  Sections~\ref{sec:zeros} and~\ref{sec:integralzero} determine the level-one zero divisor and its finite \'{e}tale lift.  Sections~\ref{sec:Legendre} and~\ref{sec:Deuring} give the Legendre and hypergeometric realization.  Section~\ref{sec:Frobenius} treats the higher Dwork approximation, elliptic precision, and direct descent to level one.  The CM classification, isogeny restrictions, and the remaining global algebraicity question are considered in Section~\ref{sec:algebraic}.

\section{Katz's geometric \texorpdfstring{$E_2$}{E2} and Frobenius}\label{sec:Katz}

We first recall the geometric interpretation of $E_2$ and
fix the normalizations used throughout the paper. Let $E_4$
and $E_6$ be the normalized Eisenstein series of weights four
and six, and put
\[
 P=E_2,\qquad Q=E_4,\qquad R=E_6.
\]
The discriminant and the modular $j$-invariant are normalized by
\[
 \Delta=\frac{E_4^3-E_6^2}{1728}=\frac{Q^3-R^2}{1728},
 \qquad
 j=\frac{E_4^3}{\Delta}=\frac{Q^3}{\Delta},
\]
We use the differential operator
\begin{equation}\label{eq:thetaDefinition}
 \thetaop=q\frac{d}{dq}.
\end{equation}
With these conventions, Ramanujan's identities take the form
\begin{align}
 \thetaop P&=\frac{P^2-Q}{12},\label{eq:RamP}\\
 \thetaop Q&=\frac{PQ-R}{3},\label{eq:RamQ}\\
 \thetaop R&=\frac{PR-Q^2}{2},\label{eq:RamR}\\
 \frac{\thetaop\Delta}{\Delta}&=P.\label{eq:RamDelta}
\end{align}
For $p\geq5$, the reduction of $E_{p-1}$ is the Hasse
invariant, whose zeros are the supersingular elliptic curves.
The congruence of Serre and Swinnerton-Dyer gives
\begin{equation}\label{eq:PtimesHasse}
 P E_{p-1}\equiv E_{p+1}\pmod p,
\end{equation}
so on the ordinary part of the curve
\begin{equation}\label{eq:Pmodpordinary}
 \overline P=\frac{\overline E_{p+1}}{\overline E_{p-1}}.
\end{equation}

Let $K_0/\Qp$ be a finite unramified extension and let $E/K_0$ have good ordinary reduction.  For a short Weierstrass equation
\begin{equation}\label{eq:Weier}
 E:\quad y^2=4x^3-g_2x-g_3,
\end{equation}
write
\begin{equation}\label{eq:omegaeta}
 \omega=\frac{dx}{y},\qquad \eta=\frac{x\,dx}{y}.
\end{equation}
The Hodge line in $H^1_{\dR}(E/K_0)$ is $K_0\omega$.  Via crystalline--de Rham comparison, the slope-zero subspace of crystalline cohomology gives a canonical line
\begin{equation}\label{eq:unitline}
 U\subset H^1_{\dR}(E/K_0)
\end{equation}
complementary to $K_0\omega$.  Hence there is a unique $u\in K_0$ such that
\begin{equation}\label{eq:ucoord}
 U=K_0(\eta+u\omega).
\end{equation}
Katz's unit-root description is given in \cite[Appendix~2]{Katz1973}; see also \cite{Katz1976}.  In the normalization used here, Mazur--Stein--Tate record the formula explicitly in \cite[\S3.2]{MST2006}: if a nonzero vector \(\gamma\omega+\delta\eta\) spans the
unit-root line, then
\[
 P(E,\omega)=-12\frac{\gamma}{\delta}.
\]
For \(U=K_0(\eta+u\omega)\), we have \(\gamma=u\) and
\(\delta=1\). Thus the formula takes the following form.

\begin{theorem}\label{thm:KatzUnit}
With the normalization above,
\begin{equation}\label{eq:KatzUnitFormula}
 P(E,\omega)=-12u.
\end{equation}
\end{theorem}

Let $\sigma$ be arithmetic Frobenius on $K_0$.  Crystalline Frobenius is $\sigma$-semilinear.  In the basis $(\omega,\eta)$ write
\begin{equation}\label{eq:Fmatrix}
 F(\omega)=\alpha\omega+\gamma\eta,
 \qquad
 F(\eta)=\beta\omega+\delta\eta.
\end{equation}
Stability of $U=K_0(\eta+u\omega)$ gives
\[
 u=\frac{\beta+\alpha\,\sigma(u)}{\delta+\gamma\,\sigma(u)}.
\]
At $u=0$, $K_0\eta$ is the slope-zero line exactly when $\beta=0$ and $\delta$ is a unit.

\begin{proposition}\label{prop:zeroUnit}
For an ordinary elliptic curve over $K_0$ written as in \eqref{eq:Weier}, the following are equivalent:
\begin{enumerate}[label=\textup{(\roman*)}]
\item $P(E,\omega)=0$;
\item $U=K_0\eta$;
\item $\eta$ spans the slope-zero Frobenius line;
\item in \eqref{eq:Fmatrix}, $\beta=0$ and $\delta\in\mathcal O_{K_0}^\times$.
\end{enumerate}
\end{proposition}

\begin{proof}
By Theorem~\ref{thm:KatzUnit}, the vanishing of $P(E,\omega)$
is equivalent to $u=0$, which is exactly $U=K_0\eta$.
The line $U$ is the slope-zero line, so this also proves
the equivalence with (iii).

If $U=K_0\eta$, Frobenius preserves this line. The coefficient
of $\omega$ in $F(\eta)$ is therefore zero, giving $\beta=0$.
Its action on the line has slope zero, so $\delta$ is a unit.
Conversely, if $\beta=0$ and $\delta$ is a unit, the line
$K_0\eta$ is Frobenius-stable and has slope zero. The
uniqueness of the slope-zero line then gives $U=K_0\eta$.
\end{proof}

The unramified hypothesis is used only for the displayed semilinear Frobenius matrix.  The criterion $E_2(E,\omega)=0$ if and only if the unit-root line is generated by $\eta$ is preserved by scalar extension and will also be used over finite extensions of $\Qp$ and over $\Cp$.

\begin{remark}\label{rem:RelativeFrobenius}
Let $\varphi$ be a Frobenius lift on a base ring and let $F:\varphi^*H\to H$ denote relative crystalline Frobenius. Write its matrix in a Hodge-adapted basis as in \eqref{eq:Fmatrix}, and put $U=\langle\eta+u\omega\rangle$. Stability of the unit-root subbundle is the identity
\[
 \beta+\alpha\varphi(u)=u\bigl(\delta+\gamma\varphi(u)\bigr).
\]
If $u(x)=0$, this identity gives
$\beta(x)=-\alpha(x)\varphi(u)(x)$. The value
$\varphi(u)(x)$ involves the parameter after applying
Frobenius, and need not be zero. Consequently,
Proposition~\ref{prop:zeroUnit} applies to the crystalline
matrix acting on the cohomology of the given curve.
A relative matrix must first be compared with that matrix.
At a point compatible with the Frobenius lift, specialization
gives the required comparison directly. At other points one
must also use the Gauss--Manin comparison between the
transformed fiber and the coefficient-conjugate fiber.
\end{remark}

\section{The ordinary zero divisor on \texorpdfstring{$X(1)$}{X(1)}}\label{sec:zeros}

Let $A=E_{p-1}$. Since $A$ is invertible on the ordinary
part of the modular curve, \eqref{eq:Pmodpordinary}
identifies the reduction of Katz's $E_2$ there with
$\overline E_{p+1}/\overline A$. We first determine the
elliptic zeros. We then describe all the other zeros by
a polynomial in $j$.

\subsection{Zeros at elliptic points}

To measure the multiplicity of a modular form at an elliptic
point, we pass to a fine modular curve of level prime to $p$
and choose a local generator of the Hodge bundle. The
multiplicity is the order of the resulting coefficient in
a local parameter on that curve. This convention is needed
because the map to the $j$-line is ramified at elliptic
points.

\begin{proposition}\label{prop:CMzeros}
Let $p\geq5$.
\begin{enumerate}[label=\textup{(\roman*)}]
\item Suppose $j(E)=0$ and $p\equiv1\pmod3$.  Then $E$ is ordinary and
\[
 P(E,\omega)=0.
\]
On a fine tame-level cover this zero has exact multiplicity two.
\item Suppose $j(E)=1728$ and $p\equiv1\pmod4$.  Then $E$ is ordinary and
\[
 P(E,\omega)=0.
\]
On a fine tame-level cover this zero is simple.
\end{enumerate}
\end{proposition}

\begin{proof}
For $j=0$ we may take $g_2=0$ in \eqref{eq:Weier}.  Over a field containing a primitive cube root $\zeta$, the automorphism
\[
 (x,y)\longmapsto(\zeta x,y)
\]
acts on de Rham cohomology by
\[
 \omega\longmapsto\zeta\omega,
 \qquad
 \eta\longmapsto\zeta^2\eta.
\]
The two eigencharacters are distinct.  The unit-root line is canonical, hence stable under every automorphism of $E$.  It cannot equal the Hodge line $K_0\omega$, so it must be $K_0\eta$.  The vanishing follows from \cref{prop:zeroUnit}.  The condition $p\equiv1\pmod3$ is exactly the ordinary condition for the $j=0$ CM curve.

To determine the order of vanishing, pass to a fine modular
curve of level prime to $p$, and then to its ordinary Igusa
tower, where the invariant differential has a canonical
trivialization.  If $q_{\mathrm{ST}}=1+t$ is a Serre--Tate coordinate centred at the point, the canonical differential trivializes the Hodge bundle and Katz's operator becomes
\[
 \thetaop=q_{\mathrm{ST}}\frac{d}{dq_{\mathrm{ST}}}
 =(1+t)\frac{d}{dt}
\]
on the resulting coefficient functions. For the normalization, see \cite[Remark~1.2.5 and Corollary~6.3.1]{HoweCircle}: we use the reciprocal of the Serre--Tate parameter in that reference, so that the sign agrees with $\thetaop=q\,d/dq$ at the cusp.  It therefore detects their exact order in the uniformizer $t$.  At $j=0$ one has $P=Q=0$ and $R\neq0$.  Ramanujan's identities give
\[
 \thetaop P=0,
 \qquad
 \thetaop^2P=\frac{R}{36}\neq0,
\]
so the zero has exact multiplicity two.

For $j=1728$ we take $g_3=0$.  The order-four automorphism
\[
 (x,y)\longmapsto(-x,iy)
\]
acts by
\[
 \omega\longmapsto i\omega,
 \qquad
 \eta\longmapsto-i\eta.
\]
Again the canonical unit-root line must be the eigenspace complementary to $K_0\omega$, namely $K_0\eta$.  The ordinary condition is $p\equiv1\pmod4$.  Here $P=R=0$ and $Q\neq0$, so
\[
 \thetaop P=-\frac{Q}{12}\neq0.
\]
Thus the zero is simple.
\end{proof}

\begin{remark}
A weight-two section vanishes at an elliptic curve whenever an automorphism acts nontrivially on the fiber of the Hodge square: invariance under the automorphism forces its value to be zero.
\end{remark}

\subsection{Reduction modulo \texorpdfstring{$p$}{p} and the zeros of \texorpdfstring{$E_{p+1}$}{E(p+1)}}

Let $A=E_{p-1}$ be the standard lift of the Hasse invariant.  On the ordinary
part of the curve $A$ is invertible, and \eqref{eq:PtimesHasse} gives
\begin{equation}\label{eq:zeroReduction}
 \overline P=\frac{\overline E_{p+1}}{\overline A}.
\end{equation}

\begin{theorem}\label{thm:Ep1zeros}
Let $p\geq5$, and let $\overline E_{p+1}$ be the reduction modulo $p$ of the
normalized Eisenstein series of weight $p+1$.
\begin{enumerate}[label=\textup{(\roman*)}]
\item Every zero of $\overline E_{p+1}$ on $X(1)_{\overline{\Fp}}$ is ordinary.
\item Every nonelliptic zero, that is, every zero with
\[
 j\notin\{0,1728\},
\]
is simple.
\item If $j=0$ is ordinary, equivalently $p\equiv1\pmod3$, then
$\overline E_{p+1}$ has exact multiplicity two at the corresponding point on a
fine tame-level cover.
\item If $j=1728$ is ordinary, equivalently $p\equiv1\pmod4$, then
$\overline E_{p+1}$ has exact multiplicity one at the corresponding point on a
fine tame-level cover.
\end{enumerate}
In particular, the usual divisor polynomial of $E_{p+1}$ modulo $p$, after the
standard elliptic factors have been removed, has only simple roots.
\end{theorem}

\begin{proof}
The reductions of the isobaric polynomials representing $E_{p-1}$ and $E_{p+1}$
are relatively prime; this is \cite[\S1.4, Theorem~5 and Corollary~1]{SerreCong}.
Katz recovers the statement by identifying the corresponding Cartier coefficients
with the reductions of $E_{p-1}$ and $E_{p+1}/12$
\cite[Theorem~3.1]{KatzZannier2021}.  Since $\overline E_{p-1}$ is the Hasse invariant and its zero
set consists precisely of the supersingular points,
$\overline E_{p+1}$ cannot vanish at a
supersingular point.  This proves (i).

Let $x$ be a nonelliptic zero.  By (i), $A(x)\neq0$, so \eqref{eq:zeroReduction}
identifies the local zero divisor of $\overline E_{p+1}$ with that of $\overline P$.
At a zero of $\overline P$, Ramanujan's first identity gives
\begin{equation}\label{eq:thetaPatzero}
 \thetaop\overline P=-\frac{Q}{12}.
\end{equation}
For $j(x)\neq0,1728$ one has $Q(x)R(x)\neq0$, while
\[
 \thetaop j=-j\frac{R}{Q}\neq0.
\]
Thus $j$ is a local parameter and
\[
 \frac{d\overline P}{dj}(x)
 =\frac{\thetaop\overline P(x)}{\thetaop j(x)}\neq0.
\]
This proves (ii).

For the elliptic points, pass to a fine tame-level cover and then to the ordinary
Igusa tower.  Write the Serre--Tate coordinate centred at the point as
$q_{\mathrm{ST}}=1+t$.  After the canonical trivialization of the Hodge bundle,
Katz's operator on coefficient functions is
\[
 \thetaop=q_{\mathrm{ST}}\frac{d}{dq_{\mathrm{ST}}}
 =(1+t)\frac{d}{dt};
\]
We use the reciprocal Serre--Tate parameter, as explained
with references in the proof of Proposition~\ref{prop:CMzeros}.
Since $1+t$ is a unit, successive applications of $\thetaop$
detect the orders one and two considered here in the
uniformizer $t$.  At $j=1728$ one has $R=0$ and $Q\neq0$, so
\eqref{eq:thetaPatzero} is nonzero and the zero is simple.  At $j=0$ one has
$Q=0$ and $R\neq0$, hence $\thetaop\overline P=0$.  Differentiating Ramanujan's
first identity once more and using $\thetaop Q=(PQ-R)/3$ gives
\[
 \thetaop^2\overline P
 =\frac{2\overline P\,\thetaop\overline P-\thetaop Q}{12}
 =\frac{R}{36}\neq0.
\]
Hence the zero at $j=0$ has exact multiplicity two.  Since $A$ is a unit at both
ordinary elliptic points, \eqref{eq:zeroReduction} transfers these multiplicities
to $\overline E_{p+1}$.  This proves (iii) and (iv).  The same elliptic exponents occur in Bianchini's independent factorization of the isobaric polynomial representing $E_{p+1}$ modulo $p$ \cite[Lemma~2 and Table~1]{Bianchini2026}.
\end{proof}

\begin{corollary}\label{cor:Ep1simple}
For every prime $p\geq5$, every zero of $\overline E_{p+1}$ with
$j\notin\{0,1728\}$ is ordinary and simple.
\end{corollary}

For the lifting argument, fix a level $N\geq3$ prime to $p$.
The corresponding modular curve $X(N)$ is a fine moduli
space, so a local family of elliptic curves is available.
Away from the elliptic points, the resulting statement
descends to the level-one curve.

\begin{theorem}\label{thm:HenselZeros}
Let $\bar x\in X(N)(\overline{\Fp})$ be a nonelliptic zero of
$\overline E_{p+1}$. Then $\bar x$ is ordinary, and its
residue disc, consisting of the geometric points reducing
to $\bar x$, contains a unique point $x$ at which $P=E_2$
vanishes. The lifted zero is simple.
\end{theorem}

\begin{proof}
Choose a finite unramified extension $K/\Qp$ whose residue
field contains $\bar x$ and its level structure. By smoothness,
the formal neighborhood of this point is described by
$\mathcal O_K[[t]]$. Choose an integral local generator
of the Hodge bundle and let $f(t)$ be the coefficient of
the Katz section with respect to its square. Then
$f\in\mathcal O_K[[t]]$, and this series converges at every
point of the residue disc $|t|<1$.

Since the Hasse invariant is a unit here, \eqref{eq:zeroReduction} identifies the reduction of $f$ with the coefficient of $\overline E_{p+1}$ multiplied by a unit. Theorem~\ref{thm:Ep1zeros} therefore gives $f(0)\in p\mathcal O_K$ and $f'(0)\in\mathcal O_K^\times$. Hensel's lemma produces $t_0\in p\mathcal O_K$ with $f(t_0)=0$.

It remains to prove uniqueness over the entire geometric
residue disc. Expand at $t_0$ and factor
\[
 f(t)=(t-t_0)V(t),\qquad V(t_0)=f'(t_0)\in\mathcal O_K^\times.
\]
The coefficients of $V$, expanded in $t-t_0$, are integral.
Its constant term is a unit. At every point of the residue
disc, $|t-t_0|<1$, so each nonconstant term has absolute value
less than one. Thus $V$ has absolute value one throughout
the disc, including at points defined over ramified extensions.
It follows that $t_0$ is the only zero and that this zero
is simple.
\end{proof}

\subsection{The supersingular polynomial and a global formula}

We now express the nonelliptic zeros in the $j$-coordinate.
The required polynomial is obtained from the divisor of
the Hasse invariant. Recall that
\[
 S_p(X)=\prod_{j_s}(X-j_s)\in\Fp[X]
\]
be the monic polynomial whose roots are the supersingular $j$-invariants in
$\overline{\Fp}$.  Write
\begin{equation}\label{eq:mdeltaepsilon}
 p-1=12m+4\delta_p+6\varepsilon_p,
 \qquad
 m=\left\lfloor\frac p{12}\right\rfloor,
\end{equation}
where
\[
 \delta_p=\begin{cases}0,&p\equiv1\pmod3,\\1,&p\equiv2\pmod3,\end{cases}
 \qquad
 \varepsilon_p=\begin{cases}0,&p\equiv1\pmod4,\\1,&p\equiv3\pmod4.\end{cases}
\]
Thus $j=0$ is supersingular exactly when $\delta_p=1$, and $j=1728$ is
supersingular exactly when $\varepsilon_p=1$.

The relation between the Hasse invariant and $S_p$ is classical.  In the notation
$Q=E_4$, $R=E_6$, write
\begin{equation}\label{eq:HasseDivisorPolynomial}
 E_{p-1}=\Delta^mQ^{\delta_p}R^{\varepsilon_p}H_p(j).
\end{equation}
Kaneko and Zagier recall that, modulo $p$ and after the monic normalization,
\begin{equation}\label{eq:KZsupersingular}
 S_p(j)=
 j^{\delta_p}(j-1728)^{\varepsilon_p}\overline{H}_p(j).
\end{equation}
In particular $S_p$ is square-free \cite{KanekoZagier1998}.

Define
\begin{equation}\label{eq:RpDefinition}
 Z_p(X)=-12S_p'(X)
 +8\delta_p\frac{S_p(X)}{X}
 +6\varepsilon_p\frac{S_p(X)}{X-1728}.
\end{equation}
Since $X\mid S_p$ when $\delta_p=1$ and $X-1728\mid S_p$ when $\varepsilon_p=1$, this expression belongs to $\Fp[X]$.

\begin{proposition}\label{prop:RpFormula}
On the ordinary part of the curve away from the elliptic points, one has in characteristic $p$
\begin{equation}\label{eq:PSpRp}
 \overline P
 =\frac{R}{Q}\,j\,\frac{Z_p(j)}{S_p(j)}.
\end{equation}
Equivalently, the nonelliptic zero divisor of $\overline P$ is the zero divisor of
$Z_p$.
\end{proposition}

\begin{proof}
Differentiate \eqref{eq:HasseDivisorPolynomial} logarithmically.  The reduction of $E_{p-1}$ has $q$-expansion $1$, hence $\thetaop\overline E_{p-1}=0$ by the $q$-expansion principle.  Ramanujan's identities give
\[
 \frac{\thetaop Q}{Q}=\frac P3-\frac{R}{3Q},
 \qquad
 \frac{\thetaop R}{R}=\frac P2-\frac{Q^2}{2R},
 \qquad
 \frac{\thetaop\Delta}{\Delta}=P,
 \qquad
 \thetaop j=-j\frac{R}{Q}.
\]
Because
\[
 m+\frac{\delta_p}{3}+\frac{\varepsilon_p}{2}
 =\frac{p-1}{12}\equiv-\frac1{12}\pmod p,
\]
the logarithmic derivative of \eqref{eq:HasseDivisorPolynomial}, followed by the
substitution \eqref{eq:KZsupersingular}, gives
\[
 P=\frac{R}{Q}
 \left(
 -12j\frac{S_p'(j)}{S_p(j)}+8\delta_p
 +6\varepsilon_p\frac{j}{j-1728}
 \right).
\]
Factoring out $j/S_p(j)$ gives \eqref{eq:PSpRp}.  Equivalently, this recovers the Kaneko--Zagier formula for $\widetilde E_{p+1}(j)$ in the normalization used here.
\end{proof}

\begin{theorem}\label{thm:RpZeros}
The polynomial $Z_p$ has the following properties.
\begin{enumerate}[label=\textup{(\roman*)}]
\item It is monic and
\[
 \deg Z_p=m+\delta_p+\varepsilon_p-1.
\]
\item It has no common root with $S_p$, and all of its roots are different from
$0$ and $1728$.
\item Every root of $Z_p$ is simple.
\item The roots of $Z_p$ are exactly the nonelliptic ordinary reduction classes
containing zeros of $P$.  Each such class contains a unique simple $p$-adic zero.
\end{enumerate}
Consequently the number of distinct ordinary residue classes
containing zeros of $P$ is
\begin{equation}\label{eq:numberZeroClasses}
 \deg Z_p+(1-\delta_p)+(1-\varepsilon_p)
 =m+1=\left\lfloor\frac p{12}\right\rfloor+1.
\end{equation}
\end{theorem}

\begin{proof}
Since $S_p$ is monic of degree $m+\delta_p+\varepsilon_p$, the coefficient of
$X^{m+\delta_p+\varepsilon_p-1}$ in \eqref{eq:RpDefinition} is
\[
 -12(m+\delta_p+\varepsilon_p)+8\delta_p+6\varepsilon_p
 =-(p-1)\equiv1\pmod p.
\]
This proves (i).

If $\alpha$ is a nonelliptic root of $S_p$, then
$Z_p(\alpha)=-12S_p'(\alpha)\neq0$ because $S_p$ is square-free.  If
$\delta_p=1$, write $S_p(X)=XT(X)$ with $T(0)\neq0$; then
$Z_p(0)=-4T(0)\neq0$.  If $\varepsilon_p=1$, write
$S_p(X)=(X-1728)T(X)$ with $T(1728)\neq0$; then
$Z_p(1728)=-6T(1728)\neq0$.  If $j=0$ is ordinary, work on a fine tame-level cover. The orders of $Q$ and $j$ there are one and three, while $R$ and $S_p(j)$ are units. Thus the prefactor $(R/Q)j/S_p(j)$ in \eqref{eq:PSpRp} has order two. By \cref{thm:Ep1zeros}, $\overline P$ has exactly this order, so $Z_p(0)\neq0$. At ordinary $j=1728$, the prefactor has order one, again equal to the exact order of $\overline P$, and therefore $Z_p(1728)\neq0$. This completes (ii).

Let $\alpha$ be a root of $Z_p$.  By (ii) it is ordinary and nonelliptic, and
\eqref{eq:PSpRp} gives $\overline P(\alpha)=0$.  Since $\overline A(\alpha)$ is a
unit, \eqref{eq:zeroReduction} gives $\overline E_{p+1}(\alpha)=0$.  By
\cref{thm:Ep1zeros}, this zero is simple.  The factors $R/Q$, $j$, and
$S_p(j)^{-1}$ in \eqref{eq:PSpRp} are all units at $\alpha$, so $Z_p$ also has a
simple root there.  This proves (iii).  Part (iv) now follows from
\cref{prop:RpFormula,thm:HenselZeros}.  Finally, the ordinary elliptic class
$j=0$ occurs exactly when $\delta_p=0$, and the class $j=1728$ occurs exactly when
$\varepsilon_p=0$, by \cref{prop:CMzeros}.  This proves \eqref{eq:numberZeroClasses}.
\end{proof}

\begin{corollary}\label{cor:exactZeroCount}
Over $\Cp$, the zeros of $P$ on elliptic curves with good
ordinary reduction have exactly
\[
 \left\lfloor\frac p{12}\right\rfloor+1
\]
distinct $j$-invariants.  All nonelliptic zeros are simple.  The elliptic zero
$j=1728$, when ordinary, is simple, while $j=0$, when ordinary, has multiplicity two
on a fine modular cover.
\end{corollary}

\begin{proof}
The assertion for the nonelliptic residue classes follows
from \cref{thm:RpZeros}. We must also exclude further zeros
in the ordinary elliptic residue classes.
By \cref{prop:CMzeros}, there is a Katz zero at the elliptic
point itself, of multiplicity two at $j=0$ and one at
$j=1728$. Theorem~\ref{thm:Ep1zeros} gives the same
multiplicities after reduction modulo $p$.

Choose an integral local parameter $t$ centered at this zero
on a fine modular cover, and let $m$ be its multiplicity.
The coefficient of the Katz section has the form $t^mV(t)$.
Indeed, the coefficients below degree $m$ vanish in
characteristic zero. The coefficient of degree $m$ is a
unit because the reduced order is also $m$. Thus $V$ has
integral coefficients and unit constant term. For $|t|<1$,
its nonconstant terms have absolute value less than one,
so $V(t)$ cannot vanish. There is therefore no additional
zero in the elliptic residue disc. Combining these classes
with the nonelliptic classes gives
\eqref{eq:numberZeroClasses}.
\end{proof}

\begin{corollary}\label{cor:criticalSS}
If $p\equiv1\pmod{12}$, then
\[
 Z_p(X)=-12S_p'(X).
\]
Hence the nonelliptic reductions of zeros of $P$ are exactly the critical points of
$S_p$.  In addition, both elliptic classes $j=0$ and $j=1728$ are ordinary zeros of
$P$.
\end{corollary}

\begin{corollary}\label{cor:fieldsZeros}
Let $\alpha$ be a root of $Z_p$ of degree $d$ over $\Fp$.  The $j$-invariant of
the unique zero of $P$ lifting the residue class of $\alpha$ is defined over the
unramified extension of $\Qp$ of degree $d$.  Thus the factorization of
$Z_p$ over $\Fp$ records the Frobenius orbits of the nonelliptic $p$-adic zeros.
\end{corollary}

\begin{proof}
Put $k=\Fp(\alpha)$ and $K=W(k)[1/p]$, where $W(k)$ is
the ring of Witt vectors of the finite field $k$. Thus $K$
is the unramified extension of $\Qp$ with residue field $k$.
Since $\alpha\neq0,1728$, the $j$-invariant is a smooth
local coordinate on the coarse modular curve. After passing
to a cover of level prime to $p$, the square of the Hodge
bundle has a local generator. The weight-zero function $E_2E_4/E_6$ descends over $W(k)$ and has the same zeros as $E_2$ in this disc, since $E_4$ and $E_6$ are units there. Its reduction has a simple zero at $\alpha$. Applying the preceding Hensel argument over $W(k)$ places the lifted invariant in $K$.

Its field over $\Qp$ is a subfield of $K$, hence unramified. Reduction of that field must contain $\alpha$, so its residue degree is at least $[k:\Fp]=d$. As $[K:\Qp]=d$, the degree is exactly $d$. Finally, applying Frobenius to a zero gives a zero above the Frobenius-conjugate residue class; uniqueness identifies it with the corresponding lift. This proves the assertion about orbits.
\end{proof}

\section{The integral zero polynomial}\label{sec:integralzero}

We collect the lifted nonelliptic invariants in a single
polynomial. The uniqueness of the lift in each residue disc
makes the set of its roots stable under Galois.

\begin{theorem}\label{thm:IntegralZeroPolynomial}
Let
\[
 d=\deg Z_p=\left\lfloor\frac p{12}\right\rfloor+\delta_p+\varepsilon_p-1.
\]
There is a monic polynomial
\begin{equation}\label{eq:IntegralZeroPolynomial}
 \mathcal Z_p(X)\in\Zp[X]
\end{equation}
of degree $d$ whose roots in $\Cp$ are exactly the nonelliptic $j$-invariants at which the Katz form $E_2$ vanishes.  It satisfies
\begin{equation}\label{eq:IntegralZeroReduction}
 \mathcal Z_p(X)\equiv Z_p(X)\pmod p,
\end{equation}
and its discriminant is a $p$-adic unit.

If
\[
 Z_p=\prod_i f_i
\]
is the factorization into monic irreducible polynomials over $\Fp$, then
\[
 \mathcal Z_p=\prod_i\widetilde f_i
\]
where $\widetilde f_i\in\Zp[X]$ is the unique monic Hensel lift of $f_i$ dividing $\mathcal Z_p$.  In particular, a factor of degree $d_i$ gives one Frobenius orbit of zeros defined over the unramified extension of $\Qp$ of degree $d_i$.
\end{theorem}

\begin{proof}
By \cref{thm:RpZeros}, every root of $Z_p$ is simple and ordinary, and each root has a unique $p$-adic lift at which $E_2$ vanishes.  Let $j_1,\ldots,j_d$ be these lifted $j$-invariants and set
\[
 \mathcal Z_p(X)=\prod_{i=1}^d(X-j_i).
\]
The set of roots is stable under $\operatorname{Gal}(\overline{\Q}_p/\Qp)$ because the Katz form and the $j$-map are defined over $\Qp$.  Each $j_i$ is integral and unramified by \cref{cor:fieldsZeros}, so the elementary symmetric functions belong to $\Zp$.  Reduction gives \eqref{eq:IntegralZeroReduction}.  Since $Z_p$ is square-free, its discriminant is nonzero modulo $p$, and therefore the discriminant of $\mathcal Z_p$ is a unit.

The final assertion is ordinary Hensel factorization.  Because the factors $f_i$ are pairwise coprime, each has a unique monic lift dividing $\mathcal Z_p$.  Frobenius acts transitively on the roots above an irreducible factor, giving the stated field of definition.
\end{proof}

Thus $\mathcal A_p=\Zp[X]/(\mathcal Z_p)$ is finite \'{e}tale with special fiber $\overline{\mathcal A}_p=\Fp[X]/(Z_p)$.

\begin{corollary}\label{cor:WittZeroAlgebra}
Compatibly with reduction, $\mathcal A_p$ is the unique finite \'{e}tale lift of $\overline{\mathcal A}_p$. If $\overline{\mathcal A}_p=\prod_i\mathbf F_{p^{d_i}}$, then
\[
 \mathcal A_p\simeq\prod_i W(\mathbf F_{p^{d_i}}).
\]
The class of $X$ is a distinguished lift of the residue $j$-coordinate, with characteristic polynomial $\mathcal Z_p$.
\end{corollary}
\begin{proof}
The factors of $Z_p$ are pairwise coprime, so Hensel
factorization and the Chinese remainder theorem identify
$\mathcal A_p$ with the product of the integer rings of
the unramified extensions of degrees $d_i$. Each such
integer ring is $W(\mathbf F_{p^{d_i}})$, the unique
unramified lift of its residue field. This proves the
assertion about the algebra and its reduction.

For $d\ge1$, in the basis $1,X,\ldots,X^{d-1}$ of
$\Zp[X]/(\mathcal Z_p)$, multiplication by the class of
$X$ has the companion matrix of $\mathcal Z_p$. Its
characteristic polynomial is therefore $\mathcal Z_p$.
For $d=0$, both algebras are zero and the characteristic
polynomial is the empty product $1$.
\end{proof}

\section{The Legendre period and differential equations}\label{sec:Legendre}

The Legendre family gives a second description of the zeros,
in which hypergeometric congruences can be applied. We first
recall the identities over $\C$, with the normalizations
needed for their reduction modulo $p$. Consider
\begin{equation}\label{eq:LegendreFamily}
 E_\lambda:\quad y^2=x(x-1)(x-\lambda),
 \qquad \lambda\in\mathbf P^1\setminus\{0,1,\infty\}.
\end{equation}
Put
\begin{equation}\label{eq:Fdef}
 F(\lambda)={}_2F_1\!\left(\frac12,\frac12;1;\lambda\right)
 =\sum_{n\geq0}\frac{(1/2)_n^2}{(n!)^2}\lambda^n.
\end{equation}
The function $F$ satisfies the Gauss equation
\begin{equation}\label{eq:PicardFuchs}
 \lambda(1-\lambda)F''+(1-2\lambda)F'-\frac14F=0.
\end{equation}
Let $\mathfrak H$ be the complex upper half-plane,
$\lambda(\tau)$ the modular lambda function, and
$\vartheta_3(\tau)=\sum_{n\in\Z}e^{\pi i n^2\tau}$ the
third Jacobi theta function. On the branch determined by
$\tau\in\mathfrak H$ and $\lambda=\lambda(\tau)$, one has
\begin{equation}\label{eq:Ftheta}
 F(\lambda(\tau))=\vartheta_3(\tau)^2.
\end{equation}
We write
\[
 C(\lambda)=F(\lambda)^2.
\]
Thus $C=\vartheta_3^4$ is a weight-two modular form on the level-two cover.

The standard theta identities give
\begin{align}
 E_4&=L(\lambda)C^2,
 \label{eq:E4Legendre}\\
 E_6&=\frac12M(\lambda)C^3,
 \label{eq:E6Legendre}\\
 \Delta&=\frac1{256}D(\lambda)C^6,
 \label{eq:DeltaLegendre}\\
 j&=256\frac{L(\lambda)^3}{D(\lambda)},
 \label{eq:jLegendre}
\end{align}
where
\begin{equation}\label{eq:LMD}
 L(\lambda)=\lambda^2-\lambda+1,
 \quad
 M(\lambda)=(\lambda+1)(\lambda-2)(2\lambda-1),
 \quad
 D(\lambda)=\lambda^2(1-\lambda)^2.
\end{equation}
The relation $E_4^3-E_6^2=1728\Delta$ becomes
\begin{equation}\label{eq:LMidentity}
 4L(\lambda)^3-M(\lambda)^2=27D(\lambda).
\end{equation}
It also follows that
\begin{equation}\label{eq:jminus1728}
 j(\lambda)-1728=64\frac{M(\lambda)^2}{D(\lambda)}.
\end{equation}

Set
\begin{equation}\label{eq:DoperatorIntro}
 \mathcal D=(1-2\lambda)+6\lambda(1-\lambda)\frac{d}{d\lambda}.
\end{equation}

Let $\thetaop=(2\pi i)^{-1}d/d\tau=q\,d/dq$.  The modular lambda function satisfies
\begin{equation}\label{eq:thetaLambda}
 \thetaop\lambda=\frac12\lambda(1-\lambda)C.
\end{equation}
Taking the logarithmic derivative of \eqref{eq:DeltaLegendre} and using
$\thetaop\Delta/\Delta=E_2$ gives the classical identity recorded by Ramanujan; see \cite[Part III, Ch.~17]{Berndt1991}.

\begin{proposition}\label{prop:E2Hypergeom}
On every Legendre period branch,
\begin{equation}\label{eq:E2FA}
 E_2=F A,
\end{equation}
where
\begin{equation}\label{eq:Adef}
 A(\lambda)
 =(1-2\lambda)F(\lambda)
 +6\lambda(1-\lambda)F'(\lambda).
\end{equation}
Equivalently,
\begin{equation}\label{eq:PhiComplex}
 \Phi(\lambda):=\frac{E_2}{C}
 =(1-2\lambda)+6\lambda(1-\lambda)\frac{F'}F.
\end{equation}
\end{proposition}

\begin{proof}
From \eqref{eq:DeltaLegendre},
\[
 \log\Delta=-\log256+2\log\lambda+2\log(1-\lambda)+6\log C.
\]
Apply $\thetaop$, use \eqref{eq:thetaLambda}, and note that $C=F^2$.  The result is
\[
 \frac{E_2}{F^2}
 =(1-2\lambda)+6\lambda(1-\lambda)\frac{F'}F,
\]
which is equivalent to \eqref{eq:E2FA}.
\end{proof}

\subsection{The Riccati equation}\label{sec:Riccati}

Set
\[
 \Phi(\lambda)=\frac{A(\lambda)}{F(\lambda)}
 =\frac{E_2}{C}.
\]
The Picard--Fuchs equation yields a first-order nonlinear equation for $\Phi$.

\begin{proposition}\label{prop:RiccatiComplex}
The function $\Phi$ satisfies
\begin{equation}\label{eq:Riccati}
 6\lambda(1-\lambda)\Phi'
 =-\Phi^2+2(1-2\lambda)\Phi-L(\lambda).
\end{equation}
\end{proposition}

\begin{proof}
Put $u=F'/F$.  Dividing \eqref{eq:PicardFuchs} by $F$ gives
\[
 \lambda(1-\lambda)(u'+u^2)+(1-2\lambda)u-\frac14=0.
\]
By \eqref{eq:PhiComplex},
\[
 u=\frac{\Phi-(1-2\lambda)}{6\lambda(1-\lambda)}.
\]
For the substitution, put $h=\lambda(1-\lambda)$ and $c=1-2\lambda$, so $h'=c$, $c'=-2$, and $\Phi=c+6hu$. The preceding equation gives
\[
 \Phi'=-2+6cu+6hu'=-\frac12-6hu^2.
\]
Multiplying by $6h$ and replacing $6hu$ by $\Phi-c$ yields
$6h\Phi'=-3h-(\Phi-c)^2$. Since $c^2+3h=L$, this is \eqref{eq:Riccati}.
\end{proof}

If $\Phi(\alpha)=0$, then
\begin{equation}\label{eq:derivativeAtZero}
 6\alpha(1-\alpha)\Phi'(\alpha)=-L(\alpha).
\end{equation}
A zero is therefore simple unless $L(\alpha)=0$, equivalently $j(\alpha)=0$.

For every nonzero solution $Y$ of the Legendre Picard--Fuchs equation, the quotient
\[
 (1-2\lambda)+6\lambda(1-\lambda)\frac{Y'}Y
\]
satisfies \eqref{eq:Riccati}. The equation is defined over $\Z[1/6]$; in characteristic $p\ge5$ it applies in particular to $Y=\Hdg_p$.

\subsection{The level-one equation}

The same differential identities give an equation directly
in the $j$-coordinate. The function
\[
 Y_p(j)=\frac{E_2^{\mathrm{Katz}}E_4}{E_6}
\]
has weight zero and is defined on the ordinary nonelliptic
part of the $j$-line. Write $\Phi_p^{\mathrm{an}}$ for the
analytic function $E_2^{\mathrm{Katz}}/C$ on the ordinary
Legendre family; its construction from hypergeometric
congruences is given in Section~\ref{sec:Frobenius}.
The pullback of $Y_p$ is
\[
 Y_p(j(T))=\frac{2L(T)}{M(T)}\Phi_p^{\mathrm{an}}(T).
\]

\begin{proposition}\label{prop:jriccati}
On this domain,
\begin{equation}\label{eq:jriccati}
 12j(j-1728)Y_p'
 =(j-1728)Y_p^2-(2j+6912)Y_p+j.
\end{equation}
In particular, at every nonelliptic zero $a$,
\begin{equation}\label{eq:jslope}
 Y_p'(a)=\frac1{12(a-1728)}\in\mathcal O_{\mathbf Q_p^{\mathrm{ur}}}^{\times}.
\end{equation}
\end{proposition}

\begin{proof}
Write $P=E_2^{\mathrm{Katz}}$, $Q=E_4$, $R=E_6$ and $Y=PQ/R$. The Ramanujan identities give
\[
 \theta Y=\frac RQ\left(-\frac{Y^2}{12}
 -\frac{Q^3}{12R^2}-\frac Y3+\frac{YQ^3}{2R^2}\right).
\]
Since $\theta j=-jR/Q$ and $Q^3/R^2=j/(j-1728)$, division and simplification give~\eqref{eq:jriccati}. Evaluation at $Y(a)=0$ gives~\eqref{eq:jslope}. The reduction of $a$ is different from $1728$, and
$p\ge5$, so the denominator in \eqref{eq:jslope} is a
unit. Here $\mathbf Q_p^{\mathrm{ur}}$ denotes the maximal
unramified extension of $\Qp$.
\end{proof}

On the residue disc of $a$, write
\[
 Y_p(j)=(j-a)V(j).
\]
In the expansion at $a$, the constant term of $V$ is $Y_p'(a)$ and all other coefficients are integral. The constant term is a unit by~\eqref{eq:jslope}, whereas every other term has positive valuation when $v_p(j-a)>0$. Thus $V$ is a unit throughout the disc. Therefore, with $v_p(p)=1$,
\begin{equation}\label{eq:distance}
 v_p(Y_p(j))=v_p(j-a).
\end{equation}
This holds for points over ramified extensions as well as unramified ones.

\section{The zero divisor on the Legendre cover}\label{sec:Deuring}

\subsection{The Deuring polynomial and the period relation}

Let $p\geq5$ be a prime and put
\[
 r=\frac{p-1}{2}.
\]
The Deuring polynomial of the Legendre family is
\begin{equation}\label{eq:Deuring}
 \Hdg_p(\lambda)
 =\sum_{k=0}^{r}\binom{r}{k}^2\lambda^k
 \in\Fp[\lambda].
\end{equation}
Its zeros are exactly the supersingular Legendre parameters, and they are simple; see, for example, \cite[Ch.~13]{Husemoller2004}.  It is also the reduction of the degree-$r$ hypergeometric truncation.

\begin{lemma}\label{lem:coefftrunc}
For $0\leq k\leq r$,
\[
 \frac{(1/2)_k^2}{(k!)^2}
 \equiv\binom{r}{k}^2\pmod p.
\]
For $r<k<p$, the same hypergeometric coefficient is divisible by $p$.
\end{lemma}

\begin{proof}
Use
\[
 \frac{(1/2)_k}{k!}=\frac1{4^k}\binom{2k}{k}.
\]
For $k\leq r$, the congruence
\[
 \frac1{4^k}\binom{2k}{k}\equiv(-1)^k\binom{r}{k}\pmod p
\]
is immediate from $r\equiv-1/2\pmod p$.  Squaring removes the sign.  If $r<k<p$, the addition $k+k$ has a carry in base $p$, so $p$ divides $\binom{2k}{k}$.
\end{proof}

\begin{proposition}\label{prop:DworkModp}
In $\Fp[[\lambda]]$ one has
\begin{equation}\label{eq:DworkModp}
 \overline F(\lambda)
 =\Hdg_p(\lambda)\overline F(\lambda^p).
\end{equation}
\end{proposition}

\begin{proof}
Write
\[
 a_n=\frac1{16^n}\binom{2n}{n}^2,
 \qquad F(\lambda)=\sum_{n\geq0}a_n\lambda^n.
\]
Let $n=ap+b$ with $0\leq b<p$.  If $b>r$, Lucas's theorem gives $a_n\equiv0\pmod p$ because the lowest base-$p$ digit in the central binomial coefficient produces a carry.  If $b\leq r$, Lucas's theorem and Fermat's congruence give
\[
 a_{ap+b}\equiv a_a a_b\pmod p.
\]
By Lemma~\ref{lem:coefftrunc}, the polynomial formed from the possible lowest digits $0\leq b\leq r$ is precisely $\Hdg_p$.  Comparing coefficients gives \eqref{eq:DworkModp}.
\end{proof}

The higher congruences will be used in
Section~\ref{sec:Frobenius}, in the form proved by Beukers
and Vlasenko \cite{BeukersVlasenko2021}. For the moment,
we apply $\mathcal D$ to the Deuring polynomial.

\begin{definition}\label{def:Bpoly}
Define
\begin{equation}\label{eq:Bpoly}
 \Bpoly_p(\lambda)
 =(1-2\lambda)\Hdg_p(\lambda)
 +6\lambda(1-\lambda)\Hdg_p'(\lambda).
\end{equation}
\end{definition}

Thus $\Bpoly_p$ is obtained from $\Hdg_p$ by the same
operator that sends $F$ to $A$. The identity
\eqref{eq:DworkModp} gives the corresponding relation
between their reductions.

\begin{proposition}\label{prop:ABreduction}
In $\Fp[[\lambda]]$,
\begin{equation}\label{eq:ABreduction}
 \overline A(\lambda)
 =\Bpoly_p(\lambda)\overline F(\lambda^p).
\end{equation}
Consequently, wherever $\Hdg_p$ is invertible,
\begin{equation}\label{eq:AFratio}
 \frac{\overline A}{\overline F}
 =\frac{\Bpoly_p}{\Hdg_p}.
\end{equation}
\end{proposition}

\begin{proof}
Apply $\mathcal D$ to \eqref{eq:DworkModp}.  In characteristic $p$ the derivative of $\overline F(\lambda^p)$ is zero.  Hence
\[
 \mathcal D\bigl(\Hdg_p(\lambda)\overline F(\lambda^p)\bigr)
 =\Bpoly_p(\lambda)\overline F(\lambda^p),
\]
which proves \eqref{eq:ABreduction}.  Divide by \eqref{eq:DworkModp} on the ordinary part of the curve.
\end{proof}

The differential equation for $\Phi$ also gives an equation
for $\Bpoly_p/\Hdg_p$ in characteristic $p$.

\begin{corollary}\label{cor:RiccatiFinite}
On the ordinary part of the curve, the rational function
\[
 \Phi_p(\lambda)=\frac{\Bpoly_p(\lambda)}{\Hdg_p(\lambda)}
\]
satisfies the reduction modulo $p$ of \eqref{eq:Riccati}:
\begin{equation}\label{eq:RiccatiFinite}
 6\lambda(1-\lambda)\Phi_p'
 =-\Phi_p^2+2(1-2\lambda)\Phi_p-L(\lambda).
\end{equation}
\end{corollary}

\begin{proof}
The Deuring polynomial satisfies the reduction modulo $p$
of the Legendre Picard--Fuchs equation. On the ordinary
part of the curve, $\Hdg_p$ is invertible, so its
logarithmic derivative is defined. We may therefore apply
the proof of Proposition~\ref{prop:RiccatiComplex} with
$F$ replaced by $\Hdg_p$. The same substitutions give
\eqref{eq:RiccatiFinite}.
\end{proof}

\subsection{The normalized Katz form}\label{sec:KatzLegendre}

On the affine Legendre curve, the section $C$ of
Section~\ref{sec:Legendre} is a nowhere-vanishing modular
form of weight two. Dividing by the appropriate power of
$C$ expresses the Hasse invariant as the Deuring polynomial.

\begin{proposition}\label{prop:HasseLegendre}
Modulo $p$,
\begin{equation}\label{eq:HasseLegendre}
 \frac{E_{p-1}}{C^r}=\Hdg_p(\lambda).
\end{equation}
\end{proposition}

\begin{proof}
Write $E_{p-1}$ as a polynomial in $E_4$ and $E_6$.  Since $p-1=2r$, every monomial $E_4^aE_6^b$ that occurs satisfies $2a+3b=r$.  Substituting
\[
 E_4=L(\lambda)C^2,\qquad E_6=\frac12M(\lambda)C^3
\]
shows that $E_{p-1}/C^r$ is a polynomial in $\lambda$ of degree at most $r$.  Its zeros are exactly the supersingular Legendre parameters.  These are simple and are also the roots of the monic degree-$r$ Deuring polynomial $\Hdg_p$.  The two polynomials therefore differ by a scalar.  Their values at the cusp $\lambda=0$ are both $1$, so the scalar is $1$.
\end{proof}

\begin{theorem}\label{thm:KatzRatio}
On the ordinary part of the Legendre curve in characteristic $p$,
\begin{equation}\label{eq:KatzRatio}
 \frac{\overline E_2}{\overline C}
 =\frac{\Bpoly_p(\lambda)}{\Hdg_p(\lambda)}.
\end{equation}
In particular, the ordinary zero divisor of Katz's $E_2$ on the affine level-two curve is the zero divisor of $\Bpoly_p$.
\end{theorem}

\begin{proof}
Identity \eqref{eq:PhiComplex} is an identity of formal $q$-expansions on the Legendre Tate curve.  The Katz $q$-expansion principle identifies its reduction with the ordinary $p$-adic Katz form.  Proposition~\ref{prop:ABreduction} then gives
\[
 \frac{\overline E_2}{\overline C}
 =\frac{\overline A}{\overline F}
 =\frac{\Bpoly_p}{\Hdg_p}.
\]
Since $C$ and $\Hdg_p$ are units on the ordinary part of the Legendre curve, the zero statement follows.
\end{proof}

\begin{proposition}\label{prop:Bbasic}
The polynomial $\Bpoly_p$ is monic of degree $(p+1)/2$.  Moreover
\[
 \gcd(\Bpoly_p,\Hdg_p)=1.
\]
Hence every root of $\Bpoly_p$ is ordinary.
\end{proposition}

\begin{proof}
The Deuring polynomial is monic of degree $r$.  The coefficient of $\lambda^{r+1}$ in \eqref{eq:Bpoly} is
\[
 -2-6r=1-3p\equiv1\pmod p,
\]
so $\Bpoly_p$ is monic of degree $r+1$.

Let $\alpha$ be a root of $\Hdg_p$.  The Deuring polynomial is square-free, and $\alpha\neq0,1$.  Thus
\[
 \Bpoly_p(\alpha)
 =6\alpha(1-\alpha)\Hdg_p'(\alpha)\neq0.
\]
This proves the coprimality.
\end{proof}

\begin{proposition}\label{prop:Bsimple}
Let $\alpha\in\Fbarp$ be a root of $\Bpoly_p$ with
\[
 \alpha\notin\{0,1\},
 \qquad L(\alpha)\neq0.
\]
Then $\alpha$ is a simple root.
\end{proposition}

\begin{proof}
By Proposition~\ref{prop:Bbasic}, $\Hdg_p(\alpha)\neq0$, so $\Phi_p=\Bpoly_p/\Hdg_p$ has a zero at $\alpha$.  Evaluating \eqref{eq:RiccatiFinite} at $\alpha$ gives
\[
 6\alpha(1-\alpha)\Phi_p'(\alpha)=-L(\alpha)\neq0.
\]
Thus $\Phi_p$, and hence $\Bpoly_p$, has a simple zero at $\alpha$.
\end{proof}

Only the roots of $L$, lying above $j=0$, can be multiple roots of $\Bpoly_p$.

\subsection{Factorization of the zero polynomial}\label{sec:factorization}

Retain the notation of Section~\ref{sec:zeros}: $S_p$ is the monic supersingular polynomial, $Z_p$ is defined by \eqref{eq:RpDefinition}, and
\[
 n=\deg S_p=m+\delta_p+\varepsilon_p,
 \qquad
 d=\deg Z_p=n-1.
\]
Away from $j=0$ and $1728$, the map from the Legendre
curve to the $j$-line is unramified. We shall compare
the two zero polynomials and determine separately the
factors at these elliptic values.

\begin{lemma}\label{lem:HassePullback}
There exists $u_p\in\Fp^\times$ such that
\begin{equation}\label{eq:HassePullback}
 \Hdg_p(\lambda)
 =u_p\,
 \frac{D(\lambda)^nS_p(j(\lambda))}
 {L(\lambda)^{2\delta_p}M(\lambda)^{\varepsilon_p}}.
\end{equation}
The right side is a polynomial.
\end{lemma}

\begin{proof}
Write the characteristic-zero Hasse lift in the standard level-one form
\[
 E_{p-1}=\Delta^mE_4^{\delta_p}E_6^{\varepsilon_p}G_p(j),
\]
where $G_p$ is a polynomial in $j$.  Divide by $C^r$ and use
\eqref{eq:E4Legendre}--\eqref{eq:DeltaLegendre}.  One obtains
\[
 \Hdg_p(\lambda)
 =256^{-m}2^{-\varepsilon_p}
 D^mL^{\delta_p}M^{\varepsilon_p}\overline G_p(j).
\]
Kaneko and Zagier's relation between $\overline G_p$ and the monic supersingular polynomial is, up to a nonzero scalar,
\[
 S_p(j)=j^{\delta_p}(j-1728)^{\varepsilon_p}\overline G_p(j).
\]
Substitute \eqref{eq:jLegendre} and \eqref{eq:jminus1728}.  The powers of $D$, $L$, and $M$ simplify to the expression in \eqref{eq:HassePullback}.  All omitted constants are nonzero modulo $p$ and are absorbed in $u_p$.
\end{proof}

We now combine this identity with the level-one formula
of Proposition~\ref{prop:RpFormula}. Away from the
supersingular and elliptic points, that formula reads
\[
 \overline E_2
 =\frac{E_6}{E_4}\,j\,\frac{Z_p(j)}{S_p(j)}.
\]
Dividing by $C$ and using \eqref{eq:E4Legendre} and \eqref{eq:E6Legendre} gives
\begin{equation}\label{eq:PhiZ}
 \frac{\overline E_2}{\overline C}
 =\frac{M(\lambda)}{2L(\lambda)}
 j(\lambda)\frac{Z_p(j(\lambda))}{S_p(j(\lambda))}.
\end{equation}

\begin{theorem}\label{thm:mainFactorization}
There exists a unique $c_p\in\Fp^\times$ such that
\begin{equation}\label{eq:mainFactorization}
 \Bpoly_p(\lambda)
 =c_p\,
 L(\lambda)^{2(1-\delta_p)}
 M(\lambda)^{1-\varepsilon_p}
 D(\lambda)^d
 Z_p(j(\lambda)).
\end{equation}
Here $d=\deg Z_p$, and $D^dZ_p(j(\lambda))$ is a polynomial of degree $6d$.
\end{theorem}

\begin{proof}
By \cref{thm:KatzRatio}, the left side divided by $\Hdg_p$ equals $\overline E_2/\overline C$.  Substitute \eqref{eq:HassePullback} and \eqref{eq:PhiZ}.  Since
\[
 j(\lambda)=256\frac{L^3}{D},
\]
all factors $S_p(j)$ cancel and the remaining powers are
\[
 D^{n-1}L^{2-2\delta_p}M^{1-\varepsilon_p}Z_p(j).
\]
Because $d=n-1$, these are the powers in \eqref{eq:mainFactorization}. To check that the rational expression is a polynomial, write $Z_p(X)=\sum_{k=0}^d z_kX^k$. Then
\[
 D^dZ_p(256L^3/D)=\sum_{k=0}^d z_k256^kL^{3k}D^{d-k}.
\]
The term $k=d$ has degree $6d$ and nonzero leading coefficient $256^d$; every other term has smaller degree. The exponents of $L$ and $M$ outside this sum are nonnegative. We therefore have an equality of polynomials, first established on the ordinary nonelliptic open set and hence everywhere. Monicity of $\Bpoly_p$ determines the scalar uniquely.
\end{proof}

\begin{corollary}\label{cor:rootStructure}
The roots of $\Bpoly_p$ have the following form.
\begin{enumerate}[label=\textup{(\roman*)}]
\item If $\alpha$ is a nonelliptic root of $Z_p$, the six points of the Legendre line above $j=\alpha$ are simple roots of $\Bpoly_p$.
\item If $j=0$ is ordinary, equivalently $\delta_p=0$, the two roots of $L$ are roots of $\Bpoly_p$ of exact multiplicity two.
\item If $j=1728$ is ordinary, equivalently $\varepsilon_p=0$, the three roots of $M$ are simple roots of $\Bpoly_p$.
\item No supersingular Legendre parameter is a root of $\Bpoly_p$.
\end{enumerate}
\end{corollary}

\begin{proof}
The final assertion is Proposition~\ref{prop:Bbasic}.  Away from $j=0,1728$, the map $\lambda\mapsto j$ has degree six and is unramified.  The nonelliptic roots of $Z_p$ are simple by the level-one zero theorem, or alternatively by combining Proposition~\ref{prop:Bsimple} with \eqref{eq:mainFactorization}.  This proves (i).  The factors $L^{2(1-\delta_p)}$ and $M^{1-\varepsilon_p}$ give the remaining multiplicities.  The factor $D^dZ_p(j(\lambda))$ is nonzero at the roots of $L$ and $M$, because $Z_p$ has no roots at the elliptic $j$-values.
\end{proof}

\section{Higher Dwork approximation and \texorpdfstring{$p$}{p}-adic lifting}\label{sec:Frobenius}

We have shown that $\overline E_2/\overline C=\Bpoly_p/\Hdg_p$
on the ordinary part of the Legendre curve. To lift this
identity modulo higher powers of $p$, we use congruences
for the logarithmic derivatives of hypergeometric
polynomials. The comparison at the cusps will also be
needed, so we specify the ring in which the congruences hold.

Write
\begin{equation}\label{eq:anDef}
 a_n=\frac{(1/2)_n^2}{(n!)^2}
 =\frac1{16^n}\binom{2n}{n}^2
\end{equation}
and, for $s\geq1$, define
\begin{equation}\label{eq:FsDef}
 F_s(\lambda)=\sum_{n=0}^{p^s-1}a_n\lambda^n\in\Zp[\lambda].
\end{equation}
Let $\widetilde\Hdg_p\in\Zp[\lambda]$ denote the coefficientwise lift of $\Hdg_p$, and let
\[
 \mathcal R_p^{\mathrm H}
 =\widehat{\Zp[\lambda,\widetilde\Hdg_p(\lambda)^{-1}]}
\]
be the $p$-adic completion. Its rigid generic fiber is the Hasse affinoid in the $\lambda$-line. The good ordinary Legendre subdomain is the generic fiber of
\[
 \widehat{\Zp[\lambda,\lambda^{-1},(1-\lambda)^{-1},
 \widetilde\Hdg_p(\lambda)^{-1}]}.
\]
It is defined by $|\lambda|=|1-\lambda|=1$ in the Hasse affinoid. The Hasse affinoid also contains the two cusp residue discs.

\begin{lemma}\label{lem:FsModp}
Put
\[
 e_s=1+p+\cdots+p^{s-1}=\frac{p^s-1}{p-1}.
\]
Then
\begin{equation}\label{eq:FsModp}
 F_s(\lambda)\equiv \Hdg_p(\lambda)^{e_s}\pmod p.
\end{equation}
In particular, every $F_s$ is a unit in $\mathcal R_p^{\mathrm H}$.
\end{lemma}

\begin{proof}
The coefficient congruence in Lemma~\ref{lem:coefftrunc}, together with Lucas's theorem, gives the digit factorization
\[
 a_{n_0+n_1p+\cdots+n_{s-1}p^{s-1}}
 \equiv a_{n_0}a_{n_1}\cdots a_{n_{s-1}}\pmod p
\]
for $0\leq n_i<p$, with $a_{n_i}\equiv0$ when $n_i>(p-1)/2$.  Hence
\[
 F_s(\lambda)
 \equiv\prod_{i=0}^{s-1}\Hdg_p(\lambda^{p^i})\pmod p.
\]
Over $\Fp$ one has $\Hdg_p(\lambda^{p^i})=\Hdg_p(\lambda)^{p^i}$, which gives \eqref{eq:FsModp}.
\end{proof}

Beukers and Vlasenko prove Dwork congruences for the period and its logarithmic derivative \cite[\S3, especially Theorem~3.2 and Remark~3.3]{BeukersVlasenko2021}.

\begin{lemma}\label{lem:BVDworkLocalization}
In the constant-term realization of \cite[Remark~3.3]{BeukersVlasenko2021}, put
\[
 q(t)=F(t^2).
\]
If $\beta_p(t)$ denotes their Hasse--Witt polynomial, then
\[
 \beta_p(t)\equiv\Hdg_p(t^2)\pmod p.
\]
Consequently there is a canonical isomorphism
\[
 \widehat{\Zp[t,\beta_p(t)^{-1}]}
 \simeq\widehat{\Zp[t,\widetilde\Hdg_p(t^2)^{-1}]}.
\]
If the ring on the right is denoted by $\mathcal R_{p,t}^{\mathrm H}$, then the involution $t\mapsto-t$ gives a direct decomposition
\begin{equation}\label{eq:evenOddDecomposition}
 \mathcal R_{p,t}^{\mathrm H}
 =\mathcal R_p^{\mathrm H}\oplus t\mathcal R_p^{\mathrm H},
 \qquad \lambda=t^2.
\end{equation}
\end{lemma}

\begin{proof}
For
\[
 g(x,y)=\frac14(x+x^{-1})(y+y^{-1}),
\]
the constant term of $g^k$ is zero for odd $k$ and, for $k=2n$, equals
\[
 \frac1{16^n}\binom{2n}{n}^2=a_n.
\]
Thus the period in \cite[Remark~3.3]{BeukersVlasenko2021} is $q(t)=F(t^2)$.  Their formula
\[
 \beta_p(t)=\sum_{k=0}^{p-1}(-1)^k\binom{p-1}{k}b_kt^k
\]
and the congruence $(-1)^k\binom{p-1}{k}\equiv1\pmod p$ give
\[
 \beta_p(t)\equiv\sum_{n=0}^{(p-1)/2}a_nt^{2n}
 \equiv\Hdg_p(t^2)\pmod p.
\]
If two elements are congruent modulo $p$, each becomes a unit after inverting the other and completing $p$-adically.  The resulting maps in both directions are inverse, which gives the stated canonical isomorphism.  Since the inverted polynomial is even, the usual decomposition
$\Zp[t]=\Zp[t^2]\oplus t\Zp[t^2]$ survives localization and completion, giving \eqref{eq:evenOddDecomposition}.
\end{proof}

\begin{proposition}\label{prop:DworkLogDerivative}
There is an analytic element
\[
 \mathcal L_p\in\mathcal R_p^{\mathrm H}
\]
whose expansion at $\lambda=0$ is $F'(\lambda)/F(\lambda)$ and such that, for every $s\geq1$,
\begin{equation}\label{eq:DworkLogDerivative}
 \mathcal L_p(\lambda)
 \equiv\frac{F_s'(\lambda)}{F_s(\lambda)}
 \pmod{p^s\mathcal R_p^{\mathrm H}}.
\end{equation}
\end{proposition}

\begin{proof}
With $q(t)=F(t^2)$, the truncation $\gamma_{2p^s}(t)$ of Beukers and Vlasenko is exactly $F_s(t^2)$.  Their logarithmic-derivative congruence following equation~(7), with $m=2p^s$, gives in $\mathcal R_{p,t}^{\mathrm H}$
\[
 \frac{q'(t)}{q(t)}
 \equiv
 \frac{d}{dt}F_s(t^2)\,F_s(t^2)^{-1}
 \pmod{p^s}.
\]
Both sides are odd under $t\mapsto-t$.  By \eqref{eq:evenOddDecomposition} they lie in $t\mathcal R_p^{\mathrm H}$. Moreover,
\[
 p^s\mathcal R_{p,t}^{\mathrm H}\cap t\mathcal R_p^{\mathrm H}
 =t\,p^s\mathcal R_p^{\mathrm H}.
\]
Thus their difference remains divisible by $p^s$ after identifying the odd summand with $\mathcal R_p^{\mathrm H}$ by $tf\mapsto f$. Since $2$ is a unit, division by $2t$ in this sense preserves the congruence, including at the cusp $t=0$; it requires no localization at $t$.  The result is
\[
 \frac{F'(\lambda)}{F(\lambda)}
 \equiv\frac{F_s'(\lambda)}{F_s(\lambda)}\pmod{p^s}
\]
as analytic elements of $\mathcal R_p^{\mathrm H}$.  This defines $\mathcal L_p$ and proves the congruence.
\end{proof}

Define
\begin{equation}\label{eq:PhiAnalytic}
 \Phi_p^{\mathrm{an}}
 =(1-2\lambda)+6\lambda(1-\lambda)\mathcal L_p
 \in\mathcal R_p^{\mathrm H}.
\end{equation}
\begin{proposition}[Comparison with the Katz section]\label{prop:KatzDworkComparison}
Let $\overline{\mathscr M}(2)$ be the compactified level-two modular stack over $\Zp$, with coarse parameter line $\mathbf P^1_\lambda$. On the finite $\lambda$-chart, its formal ordinary part has coarse model
\[
 \mathfrak U=\operatorname{Spf}\mathcal R_p^{\mathrm H}.
\]
The weight-two section $C$ trivializes the square of the Hodge bundle on this chart. The quotient $E_2^{\mathrm{Katz}}/C$ extends over $\mathfrak U$, including its cusps $\lambda=0,1$, and equals $\Phi_p^{\mathrm{an}}$. In particular, restriction to $|\lambda|=|1-\lambda|=1$ gives the normalized Katz section on the good ordinary part of the Legendre curve.
\end{proposition}

\begin{proof}
Both functions will be compared as elements of
$\mathcal R_p^{\mathrm H}$. This ring describes a domain
containing the cusp residue discs as well as the subdomain
$|\lambda|=|1-\lambda|=1$ parametrizing curves with good
ordinary reduction. The cusp discs are disjoint from that
subdomain, so the extension across them must first be checked.

Level two is not a fine moduli scheme. The statements about sections can be checked after adding full level four, which is tame since $p\geq5$, and then descended; equivalently, one works on $\overline{\mathscr M}(2)$. The residual automorphism $-1$ acts trivially on the square of the Hodge bundle. The level-two descent and the formal interpretation of ordinary modular forms are given in \cite[\S\S2.8--2.9]{Katz1973}.

On the smooth finite chart, $C$ is the integral
nowhere-vanishing weight-two section used in
Section~\ref{sec:Legendre}. It also extends at the two cusps and is nonvanishing there: in \eqref{eq:E4Legendre}, both $E_4$ and $L$ are units at $\lambda=0,1$. Thus $C$ trivializes the square of the extended Hodge bundle on the entire finite chart. After division by the appropriate power of $C$, the Hasse
invariant is $\Hdg_p$. Thus its zeros are given by
$\Hdg_p=0$, and the $p$-adic completion of the ordinary
coarse chart is precisely $\mathfrak U$.

Katz's construction identifies the unit-root expression with the integral $p$-adic form having expansion $1-24\sum_{n\geq1}\sigma_1(n)q^n$; see \cite[Appendix~A2.4, especially (A2.4.3)]{Katz1973}. Its absence of negative powers at the cusps gives a section of the extended Hodge square on the compactified ordinary part, by the definition of holomorphy at infinity and the formal interpretation in \cite[\S\S2.2, 2.9]{Katz1973}. Consequently $E_2^{\mathrm{Katz}}/C$ is an element of $\mathcal R_p^{\mathrm H}$.

At the cusp $\lambda=0$, the parameter $q^{1/2}$ and $\lambda$ differ by an integral invertible change of parameter, with $\lambda=16q^{1/2}+O(q)$. The formal identity \eqref{eq:PhiComplex} therefore identifies the expansion of this quotient with that of \eqref{eq:PhiAnalytic}. Finally the expansion map
\[
 \mathcal R_p^{\mathrm H}\longrightarrow\Zp[[\lambda]]
\]
is injective. Modulo $p$ it is the injection
$\Fp[\lambda,\Hdg_p^{-1}]\hookrightarrow\Fp[[\lambda]]$, since $\Hdg_p(0)=1$. An element with zero expansion is therefore divisible by $p$; repeating and using $p$-adic separatedness makes it zero. The two elements agree on $\mathfrak U$, and hence on the stated good ordinary subdomain.
\end{proof}

For $s\geq1$ define
\begin{equation}\label{eq:BsDef}
 \Bpoly_{p,s}(\lambda)
 =(1-2\lambda)F_s(\lambda)
 +6\lambda(1-\lambda)F_s'(\lambda)
\end{equation}
and
\begin{equation}\label{eq:PhiSDef}
 \Phi_s(\lambda)=\frac{\Bpoly_{p,s}(\lambda)}{F_s(\lambda)}.
\end{equation}

\begin{theorem}\label{thm:HigherApproximation}
For every $s\geq1$ one has in $\mathcal R_p^{\mathrm H}$
\begin{equation}\label{eq:HigherApproximation}
 \Phi_p^{\mathrm{an}}
 \equiv\Phi_s
 \pmod{p^s\mathcal R_p^{\mathrm H}}.
\end{equation}
On the ordinary part of the Legendre curve this is
\begin{equation}\label{eq:HigherE2Approximation}
 \frac{E_2^{\mathrm{Katz}}}{C}
 \equiv
 (1-2\lambda)+6\lambda(1-\lambda)\frac{F_s'}{F_s}
 \pmod{p^s}.
\end{equation}
For $s=1$, reduction modulo $p$ is \eqref{eq:KatzRatio}.
\end{theorem}

\begin{proof}
The terms $1-2\lambda$ in \eqref{eq:PhiSDef} and
\eqref{eq:PhiAnalytic} cancel. Their difference is therefore
$6\lambda(1-\lambda)(\mathcal L_p-F_s'/F_s)$.
Proposition~\ref{prop:DworkLogDerivative} makes this
difference divisible by $p^s$ in $\mathcal R_p^{\mathrm H}$.
The comparison with the Katz section gives
\eqref{eq:HigherE2Approximation}. For $s=1$, use
$F_1\equiv\Hdg_p\pmod p$ and
$\Bpoly_{p,1}\equiv\Bpoly_p\pmod p$ to obtain
\eqref{eq:KatzRatio}.
\end{proof}

\begin{lemma}\label{lem:BsModp}
For every $s\geq1$,
\begin{equation}\label{eq:BsModp}
 \Bpoly_{p,s}
 \equiv
 \Hdg_p^{e_s-1}\Bpoly_p
 \pmod p.
\end{equation}
\end{lemma}

\begin{proof}
Differentiate \eqref{eq:FsModp}.  Since $e_s\equiv1\pmod p$,
\[
 F_s'\equiv \Hdg_p^{e_s-1}\Hdg_p'\pmod p.
\]
Substitution in \eqref{eq:BsDef} gives \eqref{eq:BsModp}.
\end{proof}

The congruence of functions gives a congruence between their
zeros in every nonelliptic residue disc.

\begin{theorem}\label{thm:RootApproximation}
Let $\bar\alpha\in\Fbarp$ be a nonelliptic root of $\Bpoly_p$, and let $K_0/\Qp$ be the finite unramified extension with residue field containing $\bar\alpha$.  Then the following hold.
\begin{enumerate}[label=\textup{(\roman*)}]
\item The residue disc of $\bar\alpha$ contains a unique zero $\alpha_\infty\in\mathcal O_{K_0}$ of Katz's $p$-adic $E_2$.
\item For every $s\geq1$, the polynomial $\Bpoly_{p,s}$ has a unique zero $\alpha_s\in\mathcal O_{K_0}$ in the same residue disc.
\item The roots satisfy the congruence
\begin{equation}\label{eq:RootApproximation}
 \alpha_s\equiv\alpha_\infty\pmod{p^s}.
\end{equation}
Consequently $\alpha_s\to\alpha_\infty$ $p$-adically and
\begin{equation}\label{eq:jRootApproximation}
 j(\alpha_s)\equiv j(\alpha_\infty)\pmod{p^s}.
\end{equation}
\end{enumerate}
\end{theorem}

\begin{proof}
By Propositions~\ref{prop:Bbasic} and~\ref{prop:Bsimple}, $\Hdg_p(\bar\alpha)\neq0$ and $\Bpoly_p'(\bar\alpha)\neq0$.  The reduction of $\Phi_p^{\mathrm{an}}$ is $\Bpoly_p/\Hdg_p$, so Hensel's lemma applied to the analytic function $\Phi_p^{\mathrm{an}}$ gives a unique simple zero $\alpha_\infty$ in the residue disc.

By Lemma~\ref{lem:BsModp}, the reduction of $\Bpoly_{p,s}$ at $\bar\alpha$ has derivative
\[
 \Hdg_p(\bar\alpha)^{e_s-1}\Bpoly_p'(\bar\alpha)\neq0.
\]
Hence Hensel's lemma gives a unique root $\alpha_s$ of $\Bpoly_{p,s}$ in the same disc.

Now \cref{thm:HigherApproximation} gives
\[
 \Phi_s(\alpha_\infty)\in p^s\mathcal O_{K_0}.
\]
At the simple root $\alpha_s$, the local integral expansion factors as
\[
 \Phi_s(T)=(T-\alpha_s)V_s(T).
\]
The constant term $V_s(\alpha_s)=\Phi_s'(\alpha_s)$ is a unit, since its reduction is the nonzero derivative of $\Bpoly_p/\Hdg_p$ at $\bar\alpha$. As in the proof of Theorem~\ref{thm:HenselZeros}, $V_s$ is a unit throughout this residue disc. Evaluation at $\alpha_\infty$ therefore gives
\[
 v_p(\alpha_\infty-\alpha_s)
 =v_p\bigl(\Phi_s(\alpha_\infty)\bigr)\geq s.
\]
This proves \eqref{eq:RootApproximation}.  Since the nonelliptic hypothesis makes the map $\lambda\mapsto j$ unramified at $\bar\alpha$, its derivative is a unit there, and \eqref{eq:jRootApproximation} follows.
\end{proof}

\begin{corollary}\label{cor:RootCompatibility}
With the notation of \cref{thm:RootApproximation},
\[
 \alpha_{s+1}\equiv\alpha_s\pmod{p^s}
\]
for every $s\geq1$. Since $\Bpoly_{p,s}\in\Q[T]$, each $\alpha_s$ is algebraic over $\Q$. The congruences do not imply that $\alpha_\infty$ is algebraic over $\Q$. Its field of definition over $\Qp$ is given by the unramified lifting theorem.
\end{corollary}

\subsection{Lifting near the elliptic zeros}

At an ordinary elliptic zero, the multiplicity of the Katz
section and the ramification of the $j$-map affect the
precision of the lifted invariant. We treat the two
elliptic values separately.

\begin{theorem}\label{thm:EllipticPrecision}
Let $a$ be an ordinary Legendre parameter above an elliptic Katz zero.
\begin{enumerate}[label=\textup{(\roman*)}]
\item If $j(a)=1728$, then $\Bpoly_{p,s}$ has a unique root $a_s$ in its residue disc and
\[
 v_p(a_s-a)\geq s,\qquad v_p(j(a_s)-1728)\geq2s.
\]
\item If $j(a)=0$, then $\Bpoly_{p,s}$ has exactly two roots in its residue disc, counted with multiplicity. Each such root $b_s$ satisfies
\[
 v_p(b_s-a)\geq s/2,\qquad v_p(j(b_s))\geq3s/2.
\]
The quadratic factor containing these roots is defined over the unramified field containing $a$; its splitting field need not be unramified.
\end{enumerate}
\end{theorem}

\begin{proof}
All functions have integral local expansions. At a root of $M$, the Riccati equation gives
\[
 (\Phi_p^{\mathrm{an}})'(a)=-\frac{L(a)}{6a(1-a)},
\]
which is a unit. Consequently $\Phi_p^{\mathrm{an}}(T)=(T-a)U(T)$ with $U$ an integral unit on the residue disc. Hensel's lemma gives the unique root $a_s$ of $\Bpoly_{p,s}$ there. The function $F_s$ is a unit, so~\eqref{eq:HigherApproximation} evaluated at $a_s$ gives $v_p(\Phi_p^{\mathrm{an}}(a_s))\geq s$. Thus $v_p(a_s-a)\geq s$. The identity
\[
 j(T)-1728=64M(T)^2/D(T)
\]
and the simple roots of $M$ give the second bound.

At a root of $L$, Katz vanishing and the Riccati equation give $\Phi_p^{\mathrm{an}}(a)=(\Phi_p^{\mathrm{an}})'(a)=0$. Differentiating that equation once gives
\[
 (\Phi_p^{\mathrm{an}})''(a)=-\frac{L'(a)}{6a(1-a)}.
\]
This is a unit, as is $2$. Hence $\Phi_p^{\mathrm{an}}(T)=(T-a)^2U(T)$ with $U$ an integral unit on the disc. The reduction of $\Bpoly_{p,s}$ has exactly a double root at $\bar a$, since the other factor $\Hdg_p^{e_s-1}$ is a unit there. In a local parameter centered at $a$, Weierstrass preparation separates an invertible power series from a monic polynomial $W_s$ of degree two, whose reduction is the square of the parameter. The roots of $W_s$ all reduce to $\bar a$ and account for exactly two zeros, with multiplicity, in this disc. Equivalently, polynomial Hensel factorization selects this quadratic factor of $\Bpoly_{p,s}$. For each root $b_s$, congruence~\eqref{eq:HigherApproximation} yields
\[
 2v_p(b_s-a)=v_p(\Phi_p^{\mathrm{an}}(b_s))\geq s.
\]
Finally $j=256L^3/D$ has order three at $a$, proving $v_p(j(b_s))\geq3s/2$. The assertion about the field of the factor follows from uniqueness of Hensel factorization.
\end{proof}

\begin{remark}
The ramification indices of $j$ above $1728$ and $0$ are two and three, respectively. Above $j=0$, the selected quadratic factor is defined over the unramified field, but its individual roots may require a ramified extension.
\end{remark}

\subsection{The polynomial in the \texorpdfstring{$j$}{j}-coordinate}\label{sec:normdescent}

The nonelliptic Legendre zeros occur in groups of six,
one for each lifted $j$-invariant. We now use the
$j$-map to relate their polynomial to $\mathcal Z_p$.
Let
\[
 \Omega_p^{(0)}
 =\{\bar\lambda\in\Fbarp:\Bpoly_p(\bar\lambda)=0,
       \ j(\bar\lambda)\notin\{0,1728\}\},
\]
and let
\[
 \overline G_p(T)=\prod_{\bar\lambda\in\Omega_p^{(0)}}(T-\bar\lambda)
 \in\Fp[T].
\]
By \cref{cor:rootStructure}, this is a monic polynomial of degree $6d$ and is the monic associate of $D(T)^dZ_p(j(T))$.

\begin{lemma}\label{lem:ordinaryDworkFactor}
For every $s\geq1$, there is a unique monic factor
\[
 G_{p,s}(T)\in\Zp[T]
\]
of $\Bpoly_{p,s}(T)$ whose reduction modulo $p$ is $\overline G_p(T)$.  Its roots are the unique roots $\alpha_{s,\bar\lambda}$
in the residue discs indexed by $\Omega_p^{(0)}$.
In particular, $G_{p,s}$ has degree $6d$. It is the
factor determined by Hensel's lemma and
$\overline G_p$, independently of an ordering of the roots.
\end{lemma}

\begin{proof}
The polynomial $\prod_{\bar\lambda\in\Omega_p^{(0)}}(T-\alpha_{s,\bar\lambda})$ belongs to $\Zp[T]$: each root lies in a finite unramified extension, inertia fixes it, and arithmetic Frobenius permutes the roots according to their residue classes.  This product is monic, reduces to $\overline G_p$, and divides $\Bpoly_{p,s}$.  By \cref{lem:BsModp}, the factor $\overline G_p$ is relatively prime to the complementary factor of the reduction of $\Bpoly_{p,s}$, since its roots are ordinary and nonelliptic.  Hensel factorization therefore gives uniqueness.
\end{proof}

Put $N(T)=256L(T)^3$ and retain $D(T)=T^2(1-T)^2$, so $j=N/D$.

\begin{theorem}\label{thm:direct}
For every $s\geq1$,
\begin{equation}\label{eq:direct}
 G_{p,s}(T)\equiv256^{-d}D(T)^d\mathcal Z_p\bigl(j(T)\bigr)\pmod{p^s\Zp[T]}.
\end{equation}
Consequently, the coefficients of $G_{p,s}$ determine
those of $\mathcal Z_p$ modulo $p^s$. The linear relation
between them is injective, with the inverse relation
given below.
\end{theorem}

\begin{proof}
Let $a$ run through the roots of $\mathcal Z_p$. The polynomial
\[
 256^{-1}\bigl(N(T)-aD(T)\bigr)
\]
is monic of degree six. Its reduction has six distinct roots, because the reduction of $a$ is different from $0$ and $1728$, and $j$ is \`etale there. All six roots lift uniquely in unramified residue discs. Their $j$-invariant is exactly $a$, so they are precisely the Legendre parameters of the Katz zeros above $a$.

Thus the exact monic polynomial of all nonelliptic Legendre zeros is
\[
 \mathcal G_p(T)=256^{-d}\prod_a\bigl(N(T)-aD(T)\bigr)
 =256^{-d}D(T)^d\mathcal Z_p(j(T)).
\]
The expression is a polynomial in $\Zp[T]$ of degree $6d$. In a common finite unramified splitting field, pair its roots with those of $G_{p,s}$ by their distinct reductions. Theorem~\ref{thm:RootApproximation} gives congruence modulo $p^s$ for each pair. Taking elementary symmetric functions gives $G_{p,s}\equiv\mathcal G_p\pmod{p^s}$.

For injectivity, the polynomials $N(T)^kD(T)^{d-k}$, $0\leq k\leq d$, have distinct degrees $4d+2k$ and leading coefficients $256^k$, all units. To see injectivity also over the ring $\Zp/p^s$, suppose a nonzero linear combination vanishes and choose its largest index $k$ with nonzero coefficient. Its coefficient in degree $4d+2k$ is that nonzero coefficient times the unit $256^k$; lower-index terms cannot contribute. This is a contradiction. The same descending argument gives the inverse coefficient calculation.
\end{proof}

Write
\[
 Q_s(X)=\sum_{k=0}^dq_{k,s}X^k,\qquad q_{d,s}=1,
\]
for $\mathcal Z_p$ modulo $p^s$. The inverse relation
can be expressed by successive remainders in
$(\Zp/p^s)[T]$. With $R_d=G_{p,s}\bmod p^s$, the
coefficients and remainders satisfy, for $k=d,d-1,\ldots,0$,
\begin{align*}
 q_{k,s}&=256^{d-k}[T^{4d+2k}]R_k,\\
 R_{k-1}&=R_k-256^{-d}q_{k,s}N^kD^{d-k}.
\end{align*}
The final remainder $R_{-1}$ is zero by Theorem~\ref{thm:direct}.

\section{Algebraic zeros and isogenies}\label{sec:algebraic}

We now turn to the arithmetic relations between the zeros.
The classification of CM zeros will be used both for the
$p$-power isogeny theorem and for the bounds at arbitrary
degree. Throughout this section, $\Qbar$ denotes the
algebraic closure of $\Q$ inside $\Cp$ determined by the
chosen embedding.

\subsection{The CM case}

For a CM elliptic curve, the unit-root condition can be
compared with the two complex eigenspaces in de Rham
cohomology. This reduces the question to the following
fact about the almost holomorphic Eisenstein series.

\begin{lemma}\label{lem:almost}
The zeros in the upper half-plane of
\[
 E_2^*(\tau)=E_2(\tau)-\frac3{\pi\operatorname{Im}\tau}
\]
are precisely the $\mathrm{SL}_2(\mathbf Z)$-orbits of $i$ and $\rho=e^{2\pi i/3}$.
\end{lemma}

\begin{proof}
Since $E_2^*$ transforms with weight two, it is enough to
determine its zeros in the standard closed fundamental
domain. Write $\tau=x+iy$, where $|x|\leq1/2$ and
$|\tau|\geq1$, so $y\geq\sqrt3/2$. Put
$q_0=e^{-2\pi y}<1/200$. For $0<x<1/2$, use $|\sin(n\theta)|\leq n\sin\theta$ for $0<\theta<\pi$ and $\sigma_1(n)\leq n^2$ to obtain
\[
 \left|\sum_{n\geq2}\sigma_1(n)q_0^n\sin(2\pi nx)\right|
 \leq\sin(2\pi x)\sum_{n\geq2}n^3q_0^n
 <q_0\sin(2\pi x).
\]
For the last inequality, divide the sum by $q_0$ and
use the identity
\[
 \sum_{n\geq1}n^3q^n=\frac{q(1+4q+q^2)}{(1-q)^4},\qquad 0<q\leq1/200.
\]
Thus $\operatorname{Im}E_2(x+iy)<0$ for $0<x<1/2$; reflection gives the opposite sign for $-1/2<x<0$. Subtracting $3/(\pi y)$ does not change the imaginary part. A zero must therefore lie on $x=0$ or $x=\pm1/2$.

For $x=0$, the permitted range is $y\geq1$, and
\[
 \frac d{dy}E_2^*(iy)=48\pi\sum_{n\geq1}n\sigma_1(n)e^{-2\pi ny}
 +\frac3{\pi y^2}>0.
\]
The modular transformation fixing $i$ forces $E_2^*(i)=0$, so there is no other zero on this ray.

For $x=1/2$, put $y_0=\sqrt3/2$. Each function $y^2e^{-2\pi ny}$ decreases for $y\geq y_0$. Hence
\begin{align*}
 y^2\frac d{dy}E_2^*(1/2+iy)
 &\geq\frac3\pi-48\pi y^2\sum_{n\geq1}n\sigma_1(n)e^{-2\pi ny}\\
 &\geq\frac3\pi-36\pi\sum_{n\geq1}n^3e^{-\pi\sqrt3 n}
 >\frac3\pi-\frac\pi5>0.
\end{align*}
Indeed, the rational expression for
$\sum_{n\ge1}n^3q^n$, evaluated at $q=1/200$, is less
than $1/180$. The derivative is therefore positive on
the whole ray. The order-three stabilizer forces
$E_2^*(\rho+1)=0$, where $\rho+1=1/2+iy_0$.
This is its only zero on that ray. Translation gives
the corresponding statement on $x=-1/2$, completing
the proof.
\end{proof}

\begin{theorem}\label{thm:CMGlobalAlgebraicZeros}
Let $E/\Qbar$ have CM and good ordinary reduction at the chosen place above $p\geq5$. Then Katz's $E_2$ vanishes on $E$ if and only if $j(E)=0$ or $1728$, with the corresponding ordinary condition.
\end{theorem}

\begin{proof}
Choose an algebraic differential and its normalized equation
\[
 y^2=4x^3-g_2x-g_3,\qquad \omega=dx/y,\qquad\eta=x\,dx/y.
\]
By Proposition~\ref{prop:zeroUnit}, Katz vanishing identifies the unit-root line with $\langle\eta\rangle$. The CM algebra acts through two distinct eigenspaces, one of which is the Hodge line. Functoriality of the unit-root splitting makes $\langle\eta\rangle$ the other eigenspace. The assertion that $[\eta]$ is an eigenvector is algebraic.
Indeed, choose a nonrational CM endomorphism. Its matrix
in the algebraic basis $([\omega],[\eta])$ has entries
in $\Qbar$. The coefficient of $[\omega]$ in the image
of $[\eta]$ is zero at the chosen $p$-adic embedding.
Since this coefficient is algebraic, it is zero in
$\Qbar$ and remains zero after every complex embedding. The two CM eigenspaces over $\mathbf C$ are $H^{1,0}$ and $H^{0,1}$; consequently $[\eta]$ belongs to $H^{0,1}$.

Write the complex period lattice as $\Omega_1\mathbf Z+\Omega_2\mathbf Z$, with $\tau=\Omega_2/\Omega_1$ in the upper half-plane. Let $h_1,h_2$ be the Weierstrass quasi-periods. The periods of $\eta$ are $-h_1,-h_2$, while those of a generator of $H^{0,1}$ are proportional to $\overline\Omega_1,\overline\Omega_2$. Thus
\[
 h_2=\overline\tau h_1.
\]
Using the Legendre relation $h_1\Omega_2-h_2\Omega_1=2\pi i$ gives
\[
 h_1\Omega_1=\frac\pi{\operatorname{Im}\tau}.
\]
The standard Fourier formula $h_1\Omega_1=(\pi^2/3)E_2(\tau)$ now yields $E_2^*(\tau)=0$. Lemma~\ref{lem:almost} gives $j=0$ or $1728$. The converse follows from the extra automorphisms exactly as in Proposition~\ref{prop:CMzeros}.
\end{proof}

\begin{corollary}
Every nonelliptic Katz zero over $\Cp$ has geometric endomorphism ring $\mathbf Z$. In particular it is not a canonical CM lift of its ordinary reduction.
\end{corollary}
\begin{proof}
If a characteristic-zero elliptic curve has CM, its
$j$-invariant is an algebraic singular modulus. Over $\Cp$
it is therefore isomorphic to a CM curve defined over
$\Qbar$. Rescaling the invariant differential preserves
the vanishing of $E_2$. The theorem then forces its
$j$-invariant to be $0$ or $1728$. Thus a nonelliptic
zero has no CM and has geometric endomorphism ring
$\mathbf Z$. In particular, it cannot be a canonical
CM lift of its ordinary reduction.
\end{proof}

\subsection{Isogenies of \texorpdfstring{$p$}{p}-power degree}

The fields of definition found in Corollary~\ref{cor:fieldsZeros}
are unramified. We first explain what this implies for a
cyclic isogeny of $p$-power degree between non-CM ordinary
elliptic curves.

\begin{lemma}\label{lem:UnramifiedIsogenyChain}
Let $K/\Qp$ be finite and unramified, with $p\geq5$, and let $E$ and $E'$ be non-CM elliptic curves over $K$ with good ordinary reduction. If they are joined over $\overline K$ by a cyclic isogeny of degree $p^r$, $r\geq1$, then, after reversing its direction if necessary, its kernel is the canonical subgroup of order $p^r$. In particular the target is geometrically isomorphic to the quotient obtained by $r$ successive canonical degree-$p$ isogenies.
\end{lemma}

\begin{proof}
For an ordinary elliptic curve $A/K$, the connected--\'{e}tale sequence of its $p$-divisible group gives an inertia-stable line in $A[p](\overline K)$, namely its canonical subgroup. In a basis starting with this line, the action of inertia has the form
\[
 \left.\rho_{A,p}\right|_{I_K}
 =\begin{pmatrix}\overline\chi_p&*\\0&1\end{pmatrix},
\]
where $\overline\chi_p$ is the mod-$p$ cyclotomic character. Indeed the \'{e}tale quotient is unramified, and the determinant is cyclotomic by the Weil pairing. Since $K$ is unramified, $\overline\chi_p(I_K)=\Fp^\times$. Choose an inertia element whose cyclotomic value is not $1$. It has two distinct eigenvalues, so it fixes exactly two lines. Every $G_K$-stable line must be one of them. Consequently $A$ has at most one $K$-defined subgroup of order $p$ other than the canonical subgroup.

Let $\varphi:E\to E'$ be the given isogeny. Since the curves are non-CM and isogenous, $\operatorname{Hom}_{\overline K}(E,E')$ is a free abelian group of rank one. Galois acts on this group through $\{1,-1\}$, so $\ker\varphi$ is $G_K$-stable. Each subgroup in its unique order-$p^i$ filtration is therefore defined over $K$. Quotienting gives a chain
\[
 E=E_0\longrightarrow E_1\longrightarrow\cdots\longrightarrow E_r
\]
of $K$-defined degree-$p$ isogenies, with $E_r$ geometrically isomorphic to $E'$. All intermediate curves have good ordinary reduction. Indeed, for any prime $\ell\neq p$, their rational $\ell$-adic Tate modules are isomorphic as inertia representations to that of $E$; the N\'eron--Ogg--Shafarevich criterion gives good reduction, and ordinarity is preserved by isogeny.

By Raynaud's theorem, over a complete discrete valuation ring of mixed characteristic with absolute ramification index $e<p-1$, the generic-fiber functor on finite flat commutative group schemes killed by $p$ is fully faithful \cite[Cor.~3.3.6]{Raynaud1974}. In particular, a generic-fiber isomorphism between two such group schemes extends to an isomorphism: extend it and its inverse, and use faithfulness to identify both compositions with the identity.

For a noncanonical $K$-stable line $H_K\subset A[p]_K$, let $\mathcal H$ be its finite flat schematic closure in $\mathcal A[p]$, where $\mathcal A/\mathcal O_K$ is the good reduction model. Projection gives
\[
 \mathcal H\longrightarrow\mathcal A[p]^{\mathrm{et}}.
\]
Its generic fiber is an isomorphism because $H_K$ is distinct from the canonical line. Here $e(K)=1<p-1$, so the preceding consequence makes the projection an isomorphism over $\mathcal O_K$. Hence $\mathcal H$ is \'{e}tale. Thus each isogeny in the chain has either the canonical
connected kernel or a kernel with \'{e}tale reduction.
Duality interchanges the two cases: on the ordinary special
fiber, Frobenius has connected kernel and its dual
Verschiebung has \'{e}tale kernel. The canonical kernel
is unique, and the first paragraph gives at most one
$K$-defined noncanonical kernel. Hence an intermediate
curve has at most one $K$-defined kernel of each kind.

Suppose that two consecutive isogenies
$f_i:E_i\to E_{i+1}$ and $f_{i+1}:E_{i+1}\to E_{i+2}$
had kernels of different kinds. The kernels of $f_{i+1}$ and $\widehat f_i$ would have the same type, so uniqueness makes them equal. After identifying the two targets, their composite is therefore $\widehat f_i f_i=[p]$ and has kernel $E_i[p]\simeq(\Z/p\Z)^2$. On the other hand, the kernel of this two-step portion of the original chain is a quotient of a subgroup of the cyclic group $\ker\varphi$, so it must be cyclic. This contradiction excludes every change of type.

Consequently, either every kernel in the chain is canonical
or every kernel has \'{e}tale reduction. In the second case,
reverse the chain and take the dual isogenies. All isogenies
are then canonical. Their composite has the canonical
kernel of order $p^r$. Its geometric generic fiber is cyclic
because the connected ordinary $p$-divisible group has
height one.
\end{proof}

We can now combine the lemma with the uniqueness of the
lifted zeros.

\begin{theorem}\label{thm:PPowerExclusion}
Let $E$ and $E'$ have good ordinary reduction over $\Cp$, and suppose that Katz's $E_2$ vanishes on both. If at least one of their $j$-invariants is nonelliptic, there is no cyclic isogeny
\[
 E\longrightarrow E'
\]
of degree $p^r$ for any $r\geq1$.
\end{theorem}

\begin{proof}
Suppose such an isogeny exists. The CM classification implies that both curves are non-CM. By Corollary~\ref{cor:fieldsZeros}, their $j$-invariants lie in a common finite unramified extension $K/\Qp$. They admit models with good ordinary reduction over $K$. To see this directly, their reductions avoid $0$ and $1728$, and for $a=j(E)$ the equation
\[
 y^2=x^3-\frac{3a}{a-1728}x-\frac{2a}{a-1728}
\]
has $j$-invariant $a$ and unit discriminant. Use the same construction for $E'$. Vanishing and the existence of the geometric isogeny are unchanged by these choices.

Write $\mathscr F(a)$ for the $j$-invariant of the canonical degree-$p$ quotient of a good ordinary curve with $j$-invariant $a$. This is independent of the chosen model, is compatible with base change, and satisfies
\[
 \overline{\mathscr F(a)}=\bar a^p.
\]
These are the standard properties of the canonical quotient; see \cite[Ch.~3]{Katz1973} or \cite[\S5.1]{BorgerGurney2019}. Let $\sigma$ be the Frobenius automorphism of $K/\Qp$ whose reduction is $x\mapsto x^p$. Functoriality gives $\mathscr F(\sigma a)=\sigma\mathscr F(a)$.

By Lemma~\ref{lem:UnramifiedIsogenyChain}, after interchanging the endpoints we have $j(E')=\mathscr F^r(a)$. Its reduction is $\bar a^{p^r}$. The Katz-zero set is stable under $\sigma$, and there is only one zero $j$-invariant above this nonelliptic residue class. Therefore
\[
 \mathscr F^r(a)=\sigma^r(a).
\]
We can repeat this identity because the canonical quotient
commutes with Frobenius. No assertion about vanishing on
the intermediate canonical quotients is needed. If $\mathscr F^{rn}(a)=\sigma^{rn}(a)$, functoriality gives
\[
 \mathscr F^{r(n+1)}(a)
 =\mathscr F^r(\sigma^{rn}(a))
 =\sigma^{rn}(\mathscr F^r(a))
 =\sigma^{r(n+1)}(a).
\]
Induction proves this equality for every $n\geq1$.
Let $e=[K:\Qp]$. Since $\sigma^e$ is the identity, the
case $n=e$ gives $\mathscr F^{re}(a)=a$. The canonical
quotient of degree $p^{re}$ therefore has the same
$j$-invariant as $E$ and is geometrically isomorphic to it.
Composing the quotient with such an isomorphism gives
an endomorphism of $E$ of degree $p^{re}>1$ with cyclic
kernel. But $E$ is non-CM, so every endomorphism is
multiplication by an integer. For an integer of absolute
value greater than one, the geometric kernel is a product
of two nontrivial cyclic groups and is not cyclic. This
contradiction proves the theorem.
\end{proof}

\begin{corollary}\label{cor:PPowerIsolation}
Every $p$-power isogeny class contains at most one Katz-zero $j$-invariant.
\end{corollary}
\begin{proof}
An isogeny in characteristic zero has kernel of the form $\Z/a\Z\times\Z/b\Z$, with $a\mid b$, and factors through multiplication by $a$ followed by a cyclic isogeny of degree $b/a$. For a $p$-power isogeny between distinct nonelliptic zero $j$-invariants, this cyclic degree is a positive power of $p$, contrary to the theorem. If one endpoint is CM, both are CM, and the classification leaves only $j=0$ and $1728$. These have distinct CM fields and are not isogenous to one another.
\end{proof}

\subsection{Isogenies of arbitrary degree}\label{sec:isogenies}

For a general isogeny, the change in the unit-root coordinate
is given by V\'elu's formula. We use differentials normalized
by pullback so that the resulting identity has no additional
scaling factor.

\begin{proposition}\label{prop:isogeny}
Let $\varphi:E\longrightarrow E'$ be an isogeny of degree $n$ between elliptic curves over a finite extension of $\Qp$ with good ordinary reduction. Choose differentials with $\varphi^*\omega'=\omega$, and use the corresponding normalized short Weierstrass coordinates. Put
\[
 S_{\varphi}=\sum_{0\neq Q\in\ker\varphi}x(Q).
\]
Then
\begin{equation}\label{eq:velucohom}
 \varphi^*[\eta']=n[\eta]-S_{\varphi}[\omega],
\end{equation}
and
\begin{equation}\label{eq:isogenyP}
 E_2^{\mathrm{Katz}}(E',\omega')
 =nE_2^{\mathrm{Katz}}(E,\omega)-12S_{\varphi}.
\end{equation}
Thus, if $E$ is a Katz zero, $E'$ is a Katz zero precisely when $S_{\varphi}=0$.
\end{proposition}

\begin{proof}
V\'elu's normalized formula is
\[
 x'(\varphi(P))=x(P)+\sum_{0\neq Q\in\ker\varphi}
 \bigl(x(P+Q)-x(Q)\bigr).
\]
Since $\varphi^*\omega'=\omega$, multiplication by $\omega$ turns the left side into $\varphi^*\eta'$. For each kernel point $Q$, the term $x(P+Q)\omega$ is the translate of $\eta$. Translation preserves $\omega$ and acts as the identity on $H^1_{\dR}(E)$, so each of the $n-1$ translated terms contributes $[\eta]$. The constant terms contribute $-\sum_{Q\neq0}x(Q)[\omega]$. Together with the initial $[\eta]$, these give~\eqref{eq:velucohom}. For a modern account of the normalization see~\cite[\S2.4]{Kohel1996}.

Write the unit-root lines as $U=\langle\eta+u\omega\rangle$ and $U'=\langle\eta'+u'\omega'\rangle$. Functoriality under isogeny and~\eqref{eq:velucohom} give
\[
 \varphi^*(\eta'+u'\omega')=n\eta+(u'-S_{\varphi})\omega,
 \qquad u'=nu+S_{\varphi}.
\]
The formula $E_2^{\mathrm{Katz}}=-12u$ proves~\eqref{eq:isogenyP}. An isogeny of degree divisible by $p$ still induces an isomorphism after passing to rational crystalline cohomology, so the argument also covers that case.
\end{proof}

We next fix a cyclic degree $N$. Such an isogeny is
represented by a pair $(E,H)$, where $H\subset E$ is
cyclic of order $N$ and the target is $E/H$. Two pairs
are equivalent if an isomorphism of their source curves
carries one subgroup to the other. These pairs give the
geometric points of the open coarse modular curve $Y_0(N)$.
The datum includes neither a differential nor an identification
of the target. In particular, different kernels may give
different points even when both endpoint $j$-invariants agree.

For $N>1$, let $g_N$ be the genus of the coarse modular curve $X_0(N)$, and put
\[
 c_N=\sum_{d\mid N}\varphi\bigl(\gcd(d,N/d)\bigr),\qquad
 c_N^*=\sum_{\substack{d\mid N\\d^2\neq N}}
 \varphi\bigl(\gcd(d,N/d)\bigr),
\]
where $\varphi$ denotes Euler's totient function. Thus $c_N$ is the number of geometric cusps, and $c_N^*=c_N$ unless $N$ is a square, in which case $c_N^*=c_N-\varphi(\sqrt N)$. Define
\[
 b_N=2g_N-2+c_N^*.
\]

\begin{theorem}\label{thm:IsogenyDivisor}
There is an effective divisor $\mathscr D_N$ on $X_0(N)$, defined over $\Q$ and of degree $b_N$, with the following property. Every cyclic degree-$N$ isogeny between nonelliptic Katz zeros determines a geometric point in its support. Consequently:
\begin{enumerate}[label=\textup{(\roman*)}]
\item there are at most $b_N$ equivalence classes of such directed cyclic isogenies;
\item their endpoint invariants satisfy
\[
 [\Q(j(E),j(E')):\Q]\leq b_N;
\]
\item the degrees $e,e'$ of their reductions over $\Fp$ satisfy
\[
 \operatorname{lcm}(e,e')\leq b_N.
\]
\end{enumerate}
For prime $\ell$, one has $b_\ell=2g_\ell$, so there are at most $g_\ell$ unordered pairs of nonelliptic Katz-zero $j$-invariants joined by a cyclic degree-$\ell$ isogeny.
\end{theorem}

\begin{proof}
Consider the modular unit and its logarithmic differential
\[
 U_N(\tau)=\frac{\Delta(\tau)}{\Delta(N\tau)},\qquad
 \nu_N=\frac{dU_N}{U_N}.
\]
The function $U_N$ belongs to $\Q(X_0(N))$. Indeed, on the moduli problem of cyclic isogenies, the ratio of the source and target discriminants, using differentials normalized by $\varphi^*\omega'=\omega$, is independent of the differential and is defined over $\Q$; its complex representative is $N^{-12}U_N$. The function has no zeros or poles away from the cusps, and it is nonconstant since $U_N=q^{1-N}(1+O(q))$ at infinity.

The order at a cusp represented by $a/d$, with $d\mid N$ and $\gcd(a,d)=1$, is
\[
 m_d=\operatorname{ord}_{a/d}(U_N)
 =\frac{N-d^2}{d\gcd(d,N/d)}.
\]
This is the usual cusp-order formula applied to $\eta(\tau)^{24}/\eta(N\tau)^{24}$; see \cite[\S1]{RouseWebb2015}. It also follows by using the width $N/\gcd(d^2,N)$ and the orders $1$ and $d^2/N$ of the two discriminants in the corresponding Fourier variable. There are $\varphi(\gcd(d,N/d))$ cusps with denominator $d$, so precisely $c_N^*$ cusps have $m_d\neq0$.

In a local parameter $t$ at any cusp, write $U_N=t^m v(t)$, with $v$ a unit. Then
\[
 \nu_N=m\frac{dt}{t}+\frac{v'(t)}{v(t)}\,dt.
\]
In characteristic zero this has a simple pole when $m\neq0$, and is regular when $m=0$. There are no other poles. A nonzero meromorphic differential on a smooth projective curve of genus $g_N$ has divisor of degree $2g_N-2$. Its zero divisor $\mathscr D_N=(\nu_N)_0$ is therefore effective of degree $2g_N-2+c_N^*=b_N$ and is defined over $\Q$.

Logarithmic differentiation gives
\[
 \nu_N=2\pi i\,A_N(\tau)\,d\tau,\qquad
 A_N(\tau)=E_2(\tau)-NE_2(N\tau).
\]
Take $E=\C/(\Z+\tau\Z)$ and the cyclic subgroup generated by $1/N$. The map $z\mapsto Nz$ has target $\C/(\Z+N\tau\Z)$, with normalized target differential $dz'/N$. Let $h_1$ and $h_1'$ be the first quasi-periods for the lattices $\Z+\tau\Z$ and $\Z+N\tau\Z$. For the target differential $dz'/N$, the normalized second-kind differential is $N\wp_{\Z+N\tau\Z}(z')\,dz'$. The source cycle $z:0\to1$ maps to the target path $z':0\to N$, so the integral of its pullback is $-N^2h_1'$. The same cycle has integrals $1$ and $-h_1$ against $\omega$ and $\eta$. Thus \eqref{eq:velucohom} gives $-N^2h_1'=-Nh_1-S_\varphi$. Substituting $h_1=(\pi^2/3)E_2(\tau)$ and $h_1'=(\pi^2/3)E_2(N\tau)$ yields
\[
 S_{\varphi}=\frac{\pi^2}{3}
 \bigl(N^2E_2(N\tau)-NE_2(\tau)\bigr)
 =-\frac{\pi^2N}{3}A_N(\tau).
\]
At a point with nonelliptic source, the map from the upper half-plane to the coarse $X_0(N)$ is unramified. Thus $S_{\varphi}$ and $\nu_N$ vanish at exactly the same points there, with the same local orders. The zero divisors in question are algebraic in
characteristic zero. Their equality over $\C$ therefore
persists after extension to $\Cp$. If both curves are Katz zeros, Proposition~\ref{prop:isogeny} gives $S_{\varphi}=0$, proving the asserted containment.

An effective divisor of degree $b_N$ over $\Q$ has at most $b_N$ geometric points, proving (i). If $z$ is a closed point of its support, its contribution to the degree is at least $[\kappa(z):\Q]$. A chosen geometric point over $z$ embeds $\kappa(z)$ into $\Cp$, and both endpoint invariants belong to this image because the two $j$-maps are defined over $\Q$. Thus
\[
 [\Q(j(E),j(E')):\Q]\leq[\kappa(z):\Q]\leq b_N,
\]
which proves (ii). By Corollary~\ref{cor:fieldsZeros}, the local fields generated by the two endpoint invariants are unramified of degrees $e$ and $e'$. Put $L=\Q(j(E),j(E'))$ and take the place induced by its chosen embedding into $\Cp$. Its completion is $L_v=\Qp(j(E),j(E'))$, the compositum of those two unramified local fields. Hence
\[
 \operatorname{lcm}(e,e')=[L_v:\Qp]\leq[L:\Q]\leq b_N.
\]
The middle inequality is the usual local-degree bound for a number field. This proves (iii).

If $N=\ell$ is prime, both cusps have nonzero order, so $b_\ell=2g_\ell$. A nonelliptic Katz zero is non-CM. A cyclic self-isogeny of degree greater than one is therefore impossible, and duality pairs every isogeny under consideration with a distinct point representing the reverse direction. This proves the bound on unordered pairs.
\end{proof}

The divisor $\mathscr D_N$ depends on $N$ and is independent
of the prime $p$. Every cyclic isogeny between nonelliptic
Katz zeros must give a point of this divisor; membership
in the divisor alone does not imply vanishing at either
endpoint. For $N=11,17$, or $19$, the curve has genus one.
Thus any pair under consideration has both $j$-invariants
in a common number field of degree at most two. In these
cases, the factorization of $Z_p$ excludes pairs whose
residue degrees have least common multiple greater than two.

When $b_N=0$, the divisor is empty. This gives the following
exclusions.

\begin{corollary}\label{thm:smallisogenies}
No nonelliptic Katz zero admits a cyclic isogeny to another Katz zero of any degree in
\[
 \{2,3,4,5,7,9,13,25\}.
\]
For Katz zeros with elliptic $j$-invariants, the only possibilities at these degrees are
\[
\begin{array}{c|c}
N&j(E)=j(E')\\\hline
2,5,25&1728\\
3,7&0\\
13&0\ \text{or}\ 1728\\
4,9&\text{none}.
\end{array}
\]
The elliptic classes are subject to the ordinary condition at $p$.
\end{corollary}

\begin{proof}
For each listed $N$, one has $g_N=0$ and $c_N^*=2$, hence $b_N=0$. These values follow from the standard genus formula
\[
 g_N=1+\frac{N}{12}\prod_{\ell\mid N}(1+\ell^{-1})
 -\frac{e_2(N)}4-\frac{e_3(N)}3-\frac{c_N}2,
\]
where $e_2(N)$ and $e_3(N)$ count the roots modulo $N$ of $X^2+1$ and $X^2+X+1$, respectively; see \cite[\S2]{CsirikWetherellZieve2000}. The theorem gives the nonelliptic assertion.

Since $\mathscr D_N$ is empty, $\nu_N$ has no zero on the coarse curve. The only zeros of its pullback $A_N(\tau)d\tau$ in the upper half-plane arise from ramification at elliptic points of $\Gamma_0(N)$. Such a point represents a cyclic subgroup preserved by an extra automorphism of its source. The automorphism descends to the quotient, so both endpoints have the same elliptic $j$-invariant. The congruences defining $e_2(N)$ and $e_3(N)$ give exactly the displayed cases: neither has a root modulo $4$ or $9$, only the order-two case occurs modulo $25$, and the prime-level cases are as stated. The algebraic argument is in characteristic zero and also applies when $p$ divides $N$.
\end{proof}

The existence of a finite divisor at every cyclic degree
also has a consequence for algebraicity over $\Q$.

\begin{theorem}\label{thm:algebraicpair}
If two Katz zeros have distinct $j$-invariants and are
isogenous, both invariants belong to $\Qbar$ inside $\Cp$.
Consequently, a Katz zero whose $j$-invariant is transcendental
over $\Q$ is the only Katz-zero $j$-invariant in its
isogeny class.
\end{theorem}

\begin{proof}
If either curve is CM, both are CM, and Theorem~\ref{thm:CMGlobalAlgebraicZeros} makes their $j$-invariants algebraic. Otherwise both zeros are nonelliptic. If a cyclic isogeny of degree $n>1$ joins them, Theorem~\ref{thm:IsogenyDivisor} gives the stronger bound $[\Q(j(E),j(E')):\Q]\leq b_n$.

Any isogeny between elliptic curves in characteristic zero factors as a multiplication map followed by a cyclic isogeny, after identifying the intermediate quotient by full torsion with the source. Indeed, write its kernel as $H\simeq\mathbf Z/a\mathbf Z\times\mathbf Z/b\mathbf Z$ with $a\mid b$. Its subgroup killed by $a$ has order $a^2$, so it is all of $E[a]$. The quotient $H/E[a]$ is cyclic of order $b/a$. The multiplication map $[a]$ identifies $E/E[a]$ with $E$, and the original isogeny descends to an isogeny from this quotient with kernel $H/E[a]$. If the two $j$-invariants are distinct, this residual degree exceeds one. The first part applies.
\end{proof}

\subsection{The non-CM algebraicity question}

It remains to decide whether a nonelliptic zero can itself
have algebraic $j$-invariant. For a non-CM elliptic curve
over $\Qbar$, the unit-root criterion expresses this
question as an equality between two lines in de Rham
cohomology. We describe its relation to the period
conjecture of Ancona and Fr\u{a}\c{t}il\u{a}.

\begin{remark}\label{rem:periodCondition}
Let $E/\overline{\Q}$ be non-CM with good ordinary
reduction at a place above $p\ge5$, and let $\omega$
be a nonzero invariant differential defined over
$\overline{\Q}$. Put $M=\operatorname{End}(h^1(E))$,
the endomorphism motive of the first cohomology of $E$.
For its de Rham and crystalline groups, Ancona and
Fr\u{a}\c{t}il\u{a} obtain
\[
 G_{\mathrm{dR}}(M)\simeq PGL_2,
 \qquad
 G_{\mathrm{crys}}(M)\simeq T
\]
in the ordinary case, where $T$ is the crystalline torus.
Their period space is therefore $H_M\simeq PGL_2/T$;
see \cite[Examples~5.14 and~6.17]{AnconaFratila2025}.
Over an algebraic closure, this space parametrizes ordered
pairs of distinct lines in a two-dimensional vector space.
Let $b_E$ denote the Berthelot point determined by the
crystalline--de Rham comparison. Order the two torus
characters by their crystalline slopes. The first line
of $b_E$ is then the slope-zero line transported to
de Rham cohomology.

For the normalized short Weierstrass equation determined
by $\omega$, put $\eta=x\,dx/y$. Consider the set
\[
 D_\eta=\{(L_0,L_1)\in PGL_2/T:L_0=[\eta]\}
\]
is a proper one-dimensional algebraic subvariety, and Katz's unit-root criterion gives
\begin{equation}\label{eq:periodCondition}
 E_2^{\mathrm{Katz}}(E,\omega)=0
 \quad\Longleftrightarrow\quad
 b_E\in D_\eta.
\end{equation}
Thus the exclusion of algebraic non-CM zeros would follow
from $b_E\notin D_\eta$. For this motive, the unramified
condition of \cite[Definition~7.6]{AnconaFratila2025}
holds by Remark~7.7 of that paper. Indeed, the CM
objects in the generated tannakian category are Tate
objects, whose endomorphism field is $\Q$. This condition
is a statement about the CM objects in that category.
It is distinct from the assertion about local unramified
fields in Corollary~\ref{cor:fieldsZeros}.

Conjecture~8.3 and Proposition~8.4 of the same paper would
make $b_E$ generic in $H_M$, and hence exclude
\eqref{eq:periodCondition}. The references cited here do
not prove $b_E\notin D_\eta$ unconditionally for an
ordinary non-CM curve. Related period spaces of admissible
pairs, and an Ax--Schanuel conjecture for their variation,
are studied by Howe and Klevdal
\cite{HoweKlevdal2026,HoweKlevdalIII2026}.
\end{remark}

\begin{proposition}[Conditional nonvanishing]\label{prop:ConditionalNonvanishing}
Let $E/\Qbar$ be non-CM with good ordinary reduction at the chosen place above $p\geq5$, and let $\omega$ be a nonzero algebraic invariant differential. Suppose that $M=\operatorname{End}(h^1(E))$ satisfies Conjecture~8.3 of \cite{AnconaFratila2025}. Then
\[
 E_2^{\mathrm{Katz}}(E,\omega)\neq0.
\]
\end{proposition}
\begin{proof}
The unramified condition on $M$ holds by Remark~7.7 of
\cite{AnconaFratila2025}. Under the stated conjecture,
Proposition~8.4 of that reference therefore makes the
Berthelot point $b_E$ generic in $H_M$. The two crystalline
eigenspaces give the ordered pair of distinct lines
parametrized by $H_M\simeq PGL_2/T$, with the slope-zero
line first.

The differential $\omega$ is algebraic, so the normalized
short Weierstrass model and the line $[\eta]$ are defined
over $\Qbar$. Requiring the first line of the pair to
equal $[\eta]$ defines the proper divisor $D_\eta$
of Remark~\ref{rem:periodCondition}. A generic point
cannot belong to this divisor. Thus $b_E\notin D_\eta$,
and the equivalence \eqref{eq:periodCondition} gives
$E_2^{\mathrm{Katz}}(E,\omega)\ne0$.
\end{proof}

If this conjecture holds for $\operatorname{End}(h^1(E))$ for every non-CM algebraic $E$ with good ordinary reduction at the chosen prime, then every nonelliptic Katz-zero $j$-invariant is transcendental over $\Q$. Indeed, an algebraic such invariant would give an algebraic elliptic curve and algebraic differential, contradicting the proposition; the CM alternatives have already been excluded.

\begin{question}\label{ques:globalAlgebraicity}
Can a non-CM elliptic curve over $\overline{\Q}$ with good ordinary reduction and $j\notin\{0,1728\}$ satisfy
\[
 E_2^{\mathrm{Katz}}(E,\omega)=0
\]
for an algebraic invariant differential $\omega$?
\end{question}

For $p=5,7,13$, we have $d=0$, so there are no
nonelliptic zeros and the answer is negative. For $p=11$,
there is already one nonelliptic zero. The results of
this paper do not decide whether its $j$-invariant is
algebraic over $\Q$. The fact that it is algebraic
over $\Qp$ does not settle this question.

\section*{Acknowledgment}
The first author thanks Henri Darmon for a conversation a few years ago, during which he asked about the $p$-adic counterpart of the first author's work on the zeros of complex $E_2$ and drew his attention to the paper of Mazur, Stein, and Tate~\cite{MST2006}. This conversation was the starting point of the present work.

\end{document}